\documentclass[12pt]{article}
\usepackage{amsmath,amsthm,amssymb,graphicx,color}
\DeclareMathRadical{\sqrtsign}{symbols}{"70}{largesymbols}{"70}
\newcommand{\bb}{\mathbb}

\newcommand{\half}{{\bb H}}
\newcommand{\integers}{{\bb Z}}
\newcommand{\natls}{{\bb N}}
\newcommand{\ratls}{{\bb Q}}
\newcommand{\reals}{{\bb R}}

\newlength{\figboxwidth}
\renewcommand{\bold}[1]{\medskip \noindent {\bf #1 }\nopagebreak}

\newcommand{\isom}{\cong}

\newcommand{\cross}{\times}

\newcommand{\st}{\;\: : \;\:}         

\newcommand{\zed}{\integers}

\newcommand{\Isom}{\operatorname{Isom}}

\newcommand{\vol}{\operatorname{vol}}

\newcommand{\Area}{\operatorname{Area}}

\def\@ifundefined#1#2#3%
  {\expandafter\ifx\csname#1\endcsname\relax#2\else#3\fi}

\@ifundefined{theoremstyle}{ }{
\theoremstyle{plain} 
}
\newtheorem{theorem}{Theorem}[section]

\newtheorem{proposition}[theorem]{Proposition}
\newtheorem{lemma}[theorem]{Lemma}

\newtheorem{corollary}[theorem]{Corollary}

\@ifundefined{theoremstyle}{ }{
\theoremstyle{definition} 
}
\newtheorem{definition}[theorem]{Definition}

\newcommand{\cA}{{\cal A}}

\newcommand{\cF}{{\cal F}}
\newcommand{\cG}{{\cal G}}
\newcommand{\cH}{{\cal H}}

\mathchardef\GG="321D
\newcommand{\mc}[1]{{}}
\theoremstyle{plain}

\newcommand{\Sol}{\operatorname{Sol}(m,n)}
\newcommand{\Solv}{\operatorname{Sol}}
\newcommand{\pSol}{\operatorname{Sol}(m',n')}
\newcommand{\Solnn}{\operatorname{Sol}(n,n)}
\newcommand{\pSolnn}{\operatorname{Sol}(n',n')}
\newcommand{\X}{\operatorname{X}(m,n)}
\newcommand{\pX}{\operatorname{X}(m',n')}

\newcommand{\Bilip}{\operatorname{Bilip}}

\newcommand\inv{^{-1}}

\newcommand\G{\Gamma}

\newcommand\g{\gamma}
\newcommand\Ra{\mathbb R}

\newcommand\Za{\mathbb Z}
\newcommand\Qa{\mathbb Q}

\newcommand\Ha{\mathbb H}
\mathchardef\GG="321D

\DeclareMathOperator{\QI}{QI} \DeclareMathOperator{\Sh}{Sh}
\DeclareMathOperator{\Sim}{Sim} \DeclareMathOperator{\QSim}{Qsim}
\DeclareMathOperator{\DL}{DL}

\newtheorem{defn}[theorem]{Definition}

\begin{document}
\title{Coarse differentiation of quasi-isometries I: spaces not
quasi-isometric to Cayley graphs}
\author{Alex Eskin, David Fisher and Kevin Whyte}
\date{$ $}

\maketitle

\begin{abstract}
In this paper, we prove that certain spaces are not
quasi-isometric to Cayley graphs of finitely generated groups.  In
particular, we answer a question of Woess and prove a conjecture
of Diestel and Leader by showing that certain homogeneous graphs
are not quasi-isometric to a Cayley graph of a finitely generated
group.

This paper is the first in a sequence of papers proving results
announced in \cite{EFW0}.  In particular, this paper contains many
steps in the proofs of quasi-isometric rigidity of lattices in
$\Solv$ and of the quasi-isometry classification of lamplighter
groups.  The proofs of those results are completed in \cite{EFW1}.

The method used here is based on the idea of {\em coarse
differentiation} introduced in \cite{EFW0}.
\end{abstract}

\section{Introduction and statements of rigidity results}
\label{section:rigidity}

For any group $\G$ generated by a subset $S$ one has the
associated Cayley graph, $C_{\G}(S)$.  This is the graph with
vertex set $\G$ and edges connecting any pair of elements which
differ by right multiplication by a generator.   There is a
natural $\G$ action on $C_{\G}(S)$ by left translation.   By
giving every edge length one, the Cayley graph can be made into a
(geodesic) metric space.   The distance on $\G$ viewed as the
vertices of the Cayley graph is the {\em word metric}, defined via
the norm:

 $$\|\g\|=\inf\{\text{length of a word in the generators }S
    \text{ representing } \g \text{ in } \G.\}$$

Different sets of generators give rise to different metrics and
Cayley graphs for a group but one wants these to be equivalent.
The natural notion of equivalence in this category is {\em
quasi-isometry}:

\begin{defn}
\label{defn:qi} Let $(X,d_X)$ and $(Y,d_Y)$ be metric spaces.
Given real numbers $\kappa{\geq}1$ and $C{\geq}0$,a map
$f:X{\rightarrow}Y$ is called a {\em $(\kappa,C)$-quasi-isometry}
if
\begin{enumerate} \item
$\frac{1}{\kappa}d_X(x_1,x_2)-C{\leq}d_Y(f(x_1),f(x_2)){\leq}\kappa
d_X(x_1,x_2)+C$ for all $x_1$ and $x_2$ in $X$, and, \item the $C$
neighborhood of $f(X)$ is all of $Y$.
\end{enumerate}
\end{defn}

This paper begins the proofs of results announced in \cite{EFW0}
by developing the technique of {\em coarse differentiation} first
described there.  Proofs of some of the results in \cite{EFW0} are
continued in \cite{EFW1}. Even though quasi-isometries have no
local structure and conventional derivatives do not make sense, we
essentially construct a ``coarse derivative" that models the large
scale behavior of the quasi-isometry.

A natural question which has arisen in several contexts is whether
there exist spaces not quasi-isometric to Cayley graphs. This is
uninteresting without some assumption on homogeneity on the space,
since Cayley graphs clearly have transitive isometry group. In
this paper we prove that two types of spaces are not
quasi-isometric to Cayley graphs.   The first are non-unimodular
three dimensional solvable groups which do not admit left
invariant metrics of nonpositive curvature. The second are the
Diestel-Leader graphs, homogeneous graphs first constructed in
\cite{DL} where it was conjectured that they were not
quasi-isometric to any Cayley graph.  We prove this conjecture,
thereby answering a question raised by Woess in \cite{SW,W}.

Our work is also motivated by the program initiated by Gromov to
study finitely generated groups up to quasi-isometry \cite{Gr1,
Gr2,Gr3}.  Much interesting work has been done in this direction,
see e.g. \cite{E,EF,FM1,FM2,FM3,FS,KL,MSW,P,S1,S2,Sh,W}. For a
more detailed discussion of history and motivation, see
\cite{EFW0}.

We state our results for solvable Lie groups first as it requires
less discussion:

\begin{theorem}
\label{theorem:nolattice} Let $\Sol=\Ra{\ltimes}\Ra^2$ be a
solvable Lie group where the $\Ra$ action on $\Ra^2$ is defined by
$z{\cdot}(x,y)=(e^{mz}x,e^{-nz}y)$ for for $m,n{\in}\Ra^+$ with
$m>n$. Then there is no finitely generated group $\G$
quasi-isometric to $\Sol$.
\end{theorem}

\noindent If $m>0$ and $n<0$, then $\Sol$ admits a left invariant
metric of negative curvature.  The fact that there is no finitely
generated group quasi-isometric to $G$ in this case, provided
$m{\neq}n$, is a result of Kleiner \cite{K}, see also \cite{P2}.
When $m=n$, the group $\Solnn$ contains cocompact lattices which
are (obviously) quasi-isometric to $\Solnn$.  In the sequel to
this paper we prove that any group quasi-isometric to $\Solnn$ is
virtually a lattice in $\Solnn$ \cite{EFW1}.  Many of the partial
results in this paper hold for $m{\geq}n$ and are used in that
paper as well. Note that the assumption $m \geq n$ is only to fix
orientation and that the case $m < n$ can be reduced to this one
by changing coordinates.

We also obtain the following, which is an immediate corollary of
\cite[Theorem 5.1]{FM3} and Theorem~\ref{theorem:qisol} below:
\begin{theorem}
\label{theorem:different:sol:not:qi} $\Sol$ is quasi-isometric to
$\pSol$ if and only if $m'/m = n'/n$.
\end{theorem}

Before stating the next results, we recall a definition of the
Diestel-Leader graphs, $\DL(m,n)$. In this setting,
$m,n{\in}\Za^+$ and we assume $m{\geq}n$. Let $T_1$ and $T_2$ be
regular trees of valence $m+1$ and $n+1$ respectively. Choose
orientations on the edges of $T_1$ and $T_2$ so each vertex has
$n$ (resp. $m$) edges pointing away from it. This is equivalent to
choosing ends on these trees. We can view these orientations as
defining height functions $f_1$ and $f_2$ on the trees (the
Busemann functions for the chosen ends). If one places the point
at infinity determining $f_1$ at the bottom of the page and the
point at infinity determining $f_2$ at the top of the page, then
the trees can be drawn as:

\begin{figure}[ht]
\begin{center}
\includegraphics[width=5.5in]{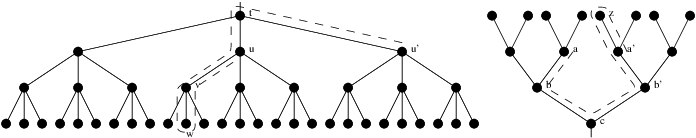}
\caption{The trees for $\DL(3,2)$. Figure borrowed from
\cite{PPS}.} \label{fig:1}
\end{center}
\end{figure}

\noindent The graph $\DL(m,n)$ is the subset of the product $T_1
\times T_2$ defined by $f_1 + f_2 = 0$. There is strong analogy
with the geometry of solvable groups which is made clear in
section $\ref{sec:solgeom}$.

\begin{theorem}
\label{theorem:dl} There is no finitely generated group
quasi-isometric to the graph $\DL(m,n)$ for $m \ne n$.
\end{theorem}

For $n=m$ the Diestel-Leader graphs arise as Cayley graphs of
lamplighter groups $\Za {\wr} F$ for $|F|=n$.  This observation
was apparently first made by R.Moeller and P.Neumann \cite{MN} and
is described explicitly, from two slightly different points of
view, in \cite{Wo2} and \cite{W}.  In \cite{EFW1} we classify
lamplighter groups up to quasi-isometry and prove that any group
quasi-isometric to a lamplighter group is a lattice {black} in
$\Isom(DL(n,n))$ for some $n$.  As discussed above, many of the
technical results in this paper are used in those proofs.

We also obtain the following analogue of
Theorem~\ref{theorem:different:sol:not:qi}:
\begin{theorem}
\label{theorem:different:dl:not:qi} If $m \neq n$ then $\DL(m,n)$ is
quasi-isometric to $\DL(m',n')$ if and only if $m$ and $m'$ are
powers of a common integer, $n$ and $n'$ are powers of a common
integer, and $\log m'/\log m = \log n'/\log n$.
\end{theorem}

\noindent  Unlike Theorem
\ref{theorem:different:sol:not:qi}, the case of this theorem where
$m=n$ is not proven in this paper. This version of the statement is
only proven in \cite{EFW1}. The case when $m=n$ here requires
additional arguments. For solvable groups, $\Solnn$ is always
quasi-isometric to $\pSolnn$ for all $n$ and $n'$.  As indicated by
the statement of the theorem, this is not true for $DL(n,n)$ and
$DL(n',n')$ which are only quasi-isometric when $n$ and $n'$ are
powers of a common integer.

The coarse differentiation approach is closely related
to results proved the method of the  ``iterated midpoint'' which is
well-known
in the theory of Banach spaces, see e.g. 
\cite{Bourgain},\cite{BL}, \cite{JLS}, \cite{Matusek}, \cite{Pr},
\cite{BJLPS}. Some results of some of those papers also have a
similar flavor, resulting in points where a map between Banach
spaces is $\epsilon$-Frechet differentiable, i.e. that the map is
sublinear distance from an affine map at some scale. The main
difference in proofs is that in our setting it is possible to
average the inequality as described in \S\ref{sec:averaging} to
obtain some control on a set of large (but not full) measure.

{\small  {\it Acknowledgements.} The first author
partially supported by NSF grant DMS-0244542, the second author was
partially supported by NSF grants DMS-0226121 and DMS-0541917 and
the third author was partially supported by NSF grant DMS-0349290
and a Sloan Foundation Fellowship.} The first two named authors
would also like to thank the Institute of Advanced Study at
Princeton and the math departments at ENS-Paris and Paris-Orsay for
hospitality while this paper was being completed. The authors also
thank Russ Lyons for useful conversations concerning Diestal-Leader
graphs and Jen Taback for comments on an early version of the
manuscript.

\section{Quasi-isometries are height respecting}
\label{section:qis}

A typical step in the study of quasi-isometric rigidity of groups
is the identification of all quasi-isometries of some space $X$
quasi-isometric to the group, see \S\ref{section:atinfinity} for
more details. For us, the space $X$ is either a solvable Lie group
$\Sol$ or $\DL(m,n)$. In all of these examples there is a special
function $h:X{\rightarrow}\Ra$ which we call the height function
and a foliation of $X$ by level sets of the height function.  We
will call a quasi-isometry of any of these spaces {\em height
respecting} if it permutes the height level sets to within bounded
distance (In \cite{FM4}, the term used is horizontal respecting).
For technical reasons, it is convenient to consider the more
general question of quasi-isometries $\Sol \to \pSol$.

For $\Sol$, the height function is $h(x,y,z)=z$.

\begin{theorem}
\label{theorem:qisol} For any $m>n>0$, any
$(\kappa,C)$-quasi-isometry $\phi: \Sol \to \pSol$ is within bounded
distance of a height respecting quasi-isometry $\hat \phi$.
Furthermore, this distance can be taken uniform in $(\kappa,C)$ and
therefore, in particular, $\hat \phi$ is a
$(\kappa',C')$-quasi-isometry where $\kappa',C'$ depend only on
$\kappa$ and $C$ and on $m$, $n$, $m'$, $n'$.
\end{theorem}

\noindent  The variant of Theorem \ref{theorem:qisol}
where $m=n$ is more difficult and is treated in \cite{EFW1}. Most of
the argument here applies in both cases and the only difference
occurs at what is labelled ``Step~II'' below. For this reason
results outside that part of this paper are all proven assuming $m
\geq n$ and not $m>n$.

In fact, Theorem \ref{theorem:qisol} can be used to identify the
self quasi-isometries of $\Sol$ completely. We will need the
following definition:

\begin{definition}[Product Map, Standard Map]
A map $\hat{\phi}: \Sol \to \pSol$ is called a {\em product map}
if it is of the form $(x,y,z) \to (f(x),g(y),q(z))$ or $(x,y,z)
\to (g(y),f(x),q(z))$, where $f$, $g$ and $q$ are functions from
$\reals \to \reals$. A product map $\hat{\phi}$ is called {\em
$b$-standard} if it is the compostion of an isometry with a map of
the form $(x,y,z) \to (f(x),g(y),z)$, where $f$ and $g$ are
Bilipshitz with the Bilipshitz constant bounded by $b$.
\end{definition}

It is known that any height-respecting quasi-isometry is at
a bounded distance from a standard map, see \cite{FM1}, and the standard maps from
$\Sol$ to $\Sol$ form a group which is isomorphic to
$(\Bilip(\reals) \cross \Bilip(\reals)){\ltimes}\Za/2\Za$ when
$m=n$ and $(\Bilip(\reals) \cross \Bilip(\reals))$ otherwise.
Given a metric space $X$, one defines $\QI(X)$ to be the group of
quasi-isometries of $X$ modulo the subgroup of those at finite
distance from the identity. Theorem~\ref{theorem:qisol} then
implies that $\QI(\Sol)=(\Bilip(\reals) \cross
\Bilip(\reals)){\rtimes}\Za/2\Za$ when $m=n$ and $(\Bilip(\reals)
\cross \Bilip(\reals))$ otherwise.  This explicit description was
conjectured by Farb and Mosher in the case $m=n$.

Recall that $\DL(m,n)$ is defined as the subset of
$T_{m+1}{\times}T_{n+1}$ where $f_1(x)+f_2(y)=0$ where $f_1$ and
$f_2$ are Busemann functions on $T_{m+1}$ and $T_{n+1}$
respectively. We fix the convention that Busemann functions
decrease as one moves toward the end from which they are defined.
We set $h((x,y))=f_m(x)=-f_n(y)$ which makes sense exactly on
$\DL(m,n){\subset}T_{m+1}{\times}T_{n+1}$.  Note that in this
choice $T_{m+1}$ branches downwards and $T_{n+1}$ branches
upwards. The reader can verify that the level sets of the height
function are orbits for a subgroup of $\Isom(DL(m,n))$.

\begin{theorem}
\label{theorem:qidl} For any $m>n$, any
$(\kappa,C)$-quasi-isometry $\varphi$ from $\DL(m,n)$ to
$\DL(m',n')$ is within bounded distance of a height respecting
quasi-isometry $\hat \varphi$. Furthermore, the bound is uniform in
$\kappa$ and $C$.
\end{theorem}

\noindent {\bf Remark:}  As above, the same result is
proven in \cite{EFW1} in the remaining case when $m=n$.

\begin{figure}
\begin{center}
\includegraphics[width=1in]{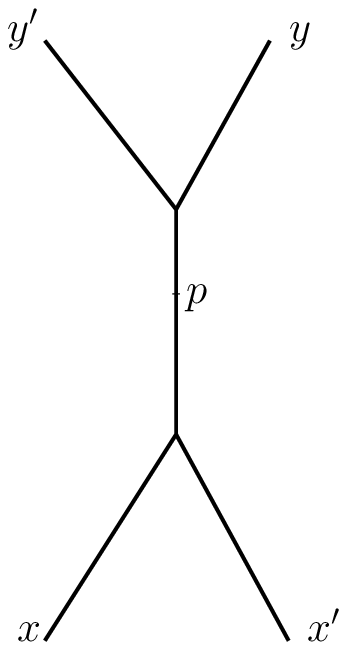}
\caption
Failure of uniqueness of the $(x,y,z)$ coordinates on $DL(m,n)$. The
point $p$ can be represented as $(x,y,z)$ or as $(x',y',z)$.

\label{fig:x}
\end{center}
\end{figure}

The discussion of standard and product maps in the setting of
$\DL(m,n)$ is slightly more complicated. We let $\Qa_l$ be the
$l$-adic rationals.  The complement of a point in the boundary at
infinity of  $T_{l+1}$ is easily seen to be  isometric to $\Qa_l$ with the $l$-adic metric.
Let $x$ be a point in $\Qa_m$ and $y$ a point in $\Qa_n$. There is a unique
vertical geodesic in $\DL(m,n)$ connecting $x$ to $y$.  To specify
a point in $\DL(m,n)$ it suffices to specify $x,y$ and a height
$z$.  We will frequently abuse notation by referring to the
$(x,y,z)$ coordinate of a point in $\DL(m,n)$ even though this
representation is highly non-unique,  see Figure~\ref{fig:x}.

Theorem \ref{theorem:qidl} can be used to identify the
quasi-isometries of $\DL(m,n)$ completely. We need to define
product and standard maps as in the case of solvable groups, but
there is an additional difficulty introduced by the non-uniqueness
of our coordinates. This is that maps of the form $(x,y,z) \to
(f(x),g(y),q(z))$, even when one assumes they are
quasi-isometries, are not well-defined. Different coordinates for
the same points will give rise to different images. We will say a
quasi-isometry $\psi$ is {\em at bounded distance} from a map of
the form $(x,y,z) \to (f(x),g(y),q(z))$ if $d(\psi(p), (f(x),
g(y), q(z)))$ is uniformly bounded for all points and all choices
$p=(x,y,z)$ of coordinates representing each point.  It is easy to
check that $(x,y,z) \to (f(x),g(y),q(z))$ is defined up to bounded
distance if we assume that the resulting map of $\DL(m,n)$ is a
quasi-isometry.  The bound depends on $\kappa,C,m,n,m'$ and $n'$.

\begin{definition}[Product Map, Standard Map]
A map $\hat{\phi}: \DL(m,n) \to \DL(m',n')$ is called a {\em
product map} if it is within bounded distance of the form $(x,y,z)
\to (f(x),g(y),q(z))$ or $(x,y,z) \to (g(y),f(x),q(z))$, where
$f: \ratls_m \to \ratls_{m'}$ (or $\ratls_{n'}$), $g: \ratls_n \to
\ratls_{n'}$ (or $\ratls_{ m'}$) and $q : \reals \to \reals$.
 A product map $\hat{\phi}$ is called {\em
$b$-standard} if it is the compostion of an isometry with a map
within bounded distance of one of the form $(x,y,z) \to
(f(x),g(y),z)$, where $f$ and $g$ are Bilipshitz with the
Bilipshitz constant bounded by $b$.
\end{definition}

\noindent Again any height-respecting quasi-isometry is at a bounded
distance from a standard map, and the standard self maps of
$\DL(m,n)$ form a group which is isomorphic to $(\Bilip(\Qa_m)
\cross \Bilip(\Qa_n)){\ltimes}\Za/2\Za$ when $m=n$ and
$(\Bilip(\Qa_m) \cross \Bilip(\Qa_n))$ otherwise.
Theorem~\ref{theorem:qisol} implies that
$\QI(\DL(m,n))=(\Bilip(\Qa_m){\times}\Bilip(\Qa_n))$ unless $m=n$
when
$\QI(\DL(m,n))=(\Bilip(\Qa_m){\times}\Bilip(\Qa_n)){\ltimes}{\Za/2\Za}$.

\section{Geometry of $\Sol$ and $\DL(m,n)$}
\label{sec:solgeom}

In this section we describe the geometry of $\Sol$ and $\DL(m,n)$,
with emphasis on the geometric facts used in our proofs.  In this
section we allow the possibility that $m=n$.   Later in
the paper we will occasionally need to develop more geometric facts
about these spaces than is described here.  We defer these facts
till later to increase readability, as they will all be isolated in
separate, clearly marked sections of the paper.

\subsection{Geodesics, quasi-geodesics and quadrilaterals}
\label{subsection:geom}

The upper half plane model of the hyperbolic plane $\Ha^2$ is the
set $\{(x, \xi) \mid \xi > 0 \}$ with the length element $ds^2 =
\frac{1}{\xi^2} (dx^2 + d\xi^2)$. If we make the change of
variable $z = \log \xi$, we get $\Ra^2$ with the length element
$ds^2 = dz^2 + e^{-2z} dx^2$. This is the {\em log model} of the
hyperbolic plane $\Ha^2$.  Note that changing $ds^2$ to
$dz^2+e^{-mz} dx^2$ we are choosing another metric of constant
negative curvature, but changing the value of the curvature. This
can be seen by checking that the substitution
$z{\rightarrow}\frac{z}{m}, x{\rightarrow}\frac{x}{m}$ is a
homothety.

The length element of $\Sol$ is:
\begin{displaymath}
ds^2 = dz^2 + e^{-2mz} dx^2 + e^{2nz} dy^2.
\end{displaymath}
Thus planes parallel to the $xz$ plane are hyperbolic planes in
the log model. Planes parallel to the $yz$ plane are {\em
upside-down} hyperbolic planes in the log model. When $m{\neq}n$,
these two families of hyperbolic planes have different
normalization on the curvature.  All of these copies of $\Ha^2$ are
isometrically embedded and totally geodesic.

\begin{itemize}
\item We use $x,y,z$ coordinates on $\Sol$, with $z$ called the {\em
    height}, and $x$ \mc{check} called the {\em depth}. The planes
  parallel to the $x z$ plane are right-side up hyperbolic planes (in
  the log model), and the planes parallel to the $y z$ plane are
  upside-down hyperbolic planes (also in the log model).

\item By ``distance'', ``area'' and ``volume''
 we mean these quantities in the $\Sol$ metric.

\end{itemize}

We will refer to lines parallel to the $x$-axis as $x$-horocycles,
and to lines parallel to the $y$-axis as $y$-horocycles. This
terminology is justified by the fact that each ($x$ or
$y$)-horocycle is indeed a horocycle in the hyperbolic plane which
contains it.





We now turn to a discussion of geodesics and quasi-geodesics in
$\Sol$.  Any geodesic in an $\Ha^2$ leaf in $\Sol$ is a geodesic.
There is a special class of geodesics, which we call {\em vertical
geodesics}.  These are the geodesics which are of the form
$\gamma(t)=(x_0,y_0, t)$ or $\gamma(t) = (x_0, y_0, -t)$. We call
the vertical geodesic  {\em upward oriented} in the first case,
and {\em downward oriented} in the second case. In both cases,
this is a unit speed parametrization. Each vertical geodesic is a
geodesic in two hyperbolic planes, the plane $y=y_0$ and the plane
$x=x_0$.

Certain quasi-geodesics in $\Sol$ are easy to describe. Given two
points $(x_0,y_0,t_0)$ and $(x_1,y_1,t_1)$, there is a geodesic
$\gamma_1$ in the hyperbolic plane $y=y_0$ that joins
$(x_0,y_0,t_0)$ to $(x_1,y_0,t_1)$ and a geodesic $\gamma_2$ in
the plane $x=x_1$ that joins $(x_1,y_0,t_1)$ to a $(x_1,y_1,t_1)$.
It is easy to check that the concatenation of $\gamma_1$ and
$\gamma_2$ is a quasi-geodesic. In first matching the $x$
coordinates and then matching the $y$ coordinates, we made a
choice.  It is possible to construct a quasi-geodesic by first
matching the $y$ coordinates and then the $x$ coordinates.  This
immediately shows that any pair of points not contained in a
hyperbolic plane in $\Sol$ can be joined by two distinct
quasi-geodesics which are not close together.  This is an aspect
of positive curvature.  One way to prove that the objects just
constructed are quasi-geodesics is to note the following: The pair
of projections $\pi_1,\pi_2:\Sol{\rightarrow}\Ha^2$ onto the $xz$
and $yz$ coordinate planes can be combined into a quasi-isometric
embedding
$\pi_1{\times}\pi_2:\Sol{\rightarrow}\Ha^2{\times}\Ha^2$.

This entire discussion is easily mimicked in $\DL(m,n)$ by
replacing geodesics and horocycles in hyperbolic planes with
geodesics and horocycles in the corresponding trees.  When we want
to state a fact that holds both for $\Sol$ and $\DL(m,n)$, we
refer to the {\em model space} which we denote by $\X$.

We define the upper boundary $\partial^+X$ as the set of
equivalence classes of vertical geodesic rays going up (where two
rays are considered equivalent if they are bounded distance
apart). The lower boundary $\partial_-X$ is defined similarly. It
is easy to see that if $X = \Sol$ case, $\partial^+X \isom \reals$
and $\partial_-X \isom \reals$. If $X = \DL(m,n)$, then
$\partial_-X \isom \ratls_m$ and $\partial^+X \isom \ratls_n$. As
discussed in Section \ref{section:qis}, if $x \in \partial_-X$, $y
\in \partial^+X$ and $z \in \reals$, we can define $(x,y,z) \in X$
as the point at height $z$ on the unique vertical geodesic
connecting $x$ and $y$.

\bold{Landau asymptotic notation.} In the following lemma and
throughout the paper, we use the notation $a = O(b)$ to mean that
$a < c_1 b$ where $c_1$ is a constant depending only on the
quasi-isometry constants $(\kappa,C)$ of $\phi$ and on the model
space or spaces (i.e. on $m$,$n$, $m'$, $n'$). We use the notation
$a = \Omega(b)$ to mean that $a > c_2 b$, where $c_2$ depends on
the same quantities as $c_1$. We also use the notation $a \GG b$
and $a \ll b$ to mean $a>C_1 b$ or $a<C_1{\inv} b$ with the same dependence
of constants.
\mc{rewrite}

We state here a key geometric fact used at various steps in the
proof.

\begin{lemma}[Quadrilaterals]
\label{lemma:quadrilateral} Let $\epsilon>0$ depending only on
$m',n'$. Suppose $p_1$, $p_2$, $q_1$, $q_2 \in \pX$ and $\gamma_{ij}:
[0,\ell_{ij}] \to \pX$ are vertical geodesic segments parametrized
by arclength.  Suppose $C > 0$ and $0 < D < \epsilon \ell_{ij}$.

Assume that for $i=1,2$, $j=1,2$,
$$d(p_i, \gamma_{ij}(0)) \leq C \qquad \text{ and  } \qquad d(q_j, \gamma_{ij}(\ell_{ij}))
\leq D,$$ so that $\gamma_{ij}$ connects the $C$-neighborhood of
$p_i$ to the $D$-neighborhood of $q_j$.  Further assume that for
$i=1,2$ and all $t$, $d(\gamma_{i1}(t),\gamma_{i2}){\geq}
(1/10) t -C$ (so that for each $i$, the two segments leaving the
neighborhood of $p_i$ diverge right away) and that for
$j=1,2$ and all $t$, $d(\gamma_{1j}(l_{1j}-t),\gamma_{2j}){\geq}
(1/10) t -D$ (so that for each $j$, the two segments leaving the
neighborhood of $q_j$ diverge right away). Then there exists $C_1 = O(C)$
and $D_1=O(D)$ such that exactly one of the following holds:
\begin{itemize}
\item[{\rm (a)}] All four $\gamma_{ij}$ are upward oriented,
$p_2$ is within $C_1$ of the $x$-horocycle passing through $p_1$,   and $q_2$ is within
$D_1$ of the $y$-horocycle passing through $q_1$.
\item[{\rm (b)}] All four $\gamma_{ij}$ are downward oriented,
$p_2$ is within $C_1$ of the $y$-horocycle passing through $p_1$,  and $q_2$ is within
$D_1$ of the $x$-horocycle passing through $q_1$.
\end{itemize}
\end{lemma}

\noindent We think of $p_1,p_2,q_1$ and $q_2$ as defining a
quadrilateral.  The content of the lemma is that any quadrilateral
has its four "corners" in pairs that lie essentially along
horocycles.

In particular, if we take a quadrilateral with geodesic segments
$\gamma_{ij}$ and with $h(p_1)=h(p_2)$ and $h(q_1)=h(q_2)$ and map
it forward under a $(\kappa,C)$-quasi-isometry $\phi:\X \to \pX$,
and if we would somehow know that $\phi$ sends each of the four
$\gamma_{ij}$ close to a vertical geodesic, then
Lemma~\ref{lemma:quadrilateral} would imply that $\phi$ sends the
$p_i$ to a pair of points at roughly the same height.

To prove Lemma~\ref{lemma:quadrilateral}, we require a
combinatorial lemma.

\begin{lemma}[Complete Bipartite Graphs]
\label{lemma:complete:bipartite:graphs} Let $\Gamma$ be an
oriented graph with four vertices $p_1, p_2, q_1, q_2$ and four
edges, such that there is exactly one edge connecting each $p_i$
to each $q_j$. Then exactly one of the following is true:
\begin{itemize}
\item[{\rm (i)}] All the edges of $\Gamma$ are from some $p_i$ to
some
  $q_j$.
\item[{\rm (ii)}] All the edges of $\Gamma$ are from some $q_j$ to
some
  $p_i$.
\item[{\rm (iii)}] There exist two vertices $v_1$ and $v_2$ which
are
  connected by two distinct directed paths.
\end{itemize}
\end{lemma}

\bold{Proof.} Since there are only 16 possibilities for $\Gamma$,
one can check directly. One way to organize the check is to let
$k$ denote the sum of number of edges outgoing from $p_1$ and the
number of edges outgoing from $p_2$. If $k=0$, (ii) holds, and if
$k=4$, then (i) holds. It is easy to check that for $1 \le k \le
3$, (iii) holds. \qed\medskip

\bold{Proof of Lemma~\ref{lemma:quadrilateral}.} Let us assume for
the moment that all the geodesics are downward oriented. Let
$x_{ij}$, $y_{ij}$ denote the $x$ and $y$ coordinates of the
vertical geodesics $\gamma_{ij}$. By the assumptions near $p_i$ we
have for $i=1,2$,
\begin{equation}
\label{eq:quadrilateralone}
C_2^{-1} \le \ln |x_{i1} - x_{i2}| e^{- m' h(p_i)}
\le C_2
\end{equation}
\noindent where $C_2 = O(C)$. The upper bound comes from the fact that $\gamma_{i1}$ and $\gamma_{i2}$
come close to $p_i$, the lower bound comes from assumption of fast divergence.
By the assumptions near $q_j$ we
have, for similar reasons, that for $j=1,2$,
\begin{equation}
\label{eq:quadrilateraltwo}
\ln |x_{1j} - x_{2j}| e^{ - m' h(q_j)} \le D_1
\end{equation}
where $D_1 = O(D)$.  Note that since for all $i,j$, $D < \epsilon
\ell_{ij}$, and so the geodesics travel a downward a long way relative to $D$,
we have
\begin{equation}
\label{eq:quadrilateralthree}
D_1 e^{m' h(q_j)} \ll C_2 e^{m' h(p_i)}.
\end{equation}

Combining the inequalities (\ref{eq:quadrilateralone}),
(\ref{eq:quadrilateraltwo}) and (\ref{eq:quadrilateralthree}), we see
that $e^{- m'(h(p_1)-h(p_2))} = O(C_2)$,
and also that $\ln |x_{1j} - x_{2j}| e^{- m' h(p_1)} = O(C_2)$. This
proves the lemma under the assumption of downward orientation.
\mc{explain more?}

The case where all the vertical geodesics are upward oriented is
identical (except that one considers differences in
$y$-coordinates instead).

To reduce to the cases already considered, we apply
Lemma~\ref{lemma:complete:bipartite:graphs} to the graph $\Gamma$
consisting of the vertices $p_1$, $p_2$, $q_1$, $q_2$ with edges
the vertical geodesics ``almost" connecting them. Suppose that
possibility (iii) of Lemma~\ref{lemma:complete:bipartite:graphs}
holds. Then, we would then have two distinct oriented paths
$\eta_1$ and $\eta_2$ connecting $v_1$ and $v_2$. Each $\eta_i$ is
either a vertical geodesics, a concatenation two vertical
geodesics, one of which ends near the beginning of the other, or a
similar concatenation of three vertical geodesics.  In each case
it is easy to check that each $\eta_i$ is close to a vertical
geodesic $\lambda_i'$. (See Lemma
\ref{lemma:epsilon:monotone:segment} for a more general variant of
this fact.) But this is a contradiction in view of the divergence
assumptions, since any pair of vertical geodesics beginning and
ending near the same point are close for their entire length. Thus
either (i) or (ii) of Lemma~\ref{lemma:complete:bipartite:graphs}
holds.

 \qed\medskip


\subsection{Volume and measure}
\label{subsection:volumeandmeasure}

There is a large difference between the unimodular and
nonunimodular examples we consider that has to do with the
measures of sets, unimodularity and amenability.  In the cases
where $m=n$ the spaces we consider are metrically amenable and
have unimodular isometry group.  When $m\neq n$ the spaces are not
metrically amenable and the isometry groups are not unimodular,
though the isometry group remains amenable as a group. \mc{Say
more or cut even this?} In particular it is immediately clear that
$\DL(n,n)$ cannot be quasi-isometric to $\DL(m,n)$ with $m \ne n$
(since one has metric F\"olner sets and the other does not). For
the same reason, $\Solnn$ is not quasi-isometric to $\pSol$ with
$m'\ne n'$.

The natural volume $\vol$ on $\DL(m,n)$ is the counting measure,
the natural volume on $\Sol$ is $\vol=e^{(n-m)z}dxdydz$. Note that
for the unimodular case where $m=n$, the volume on $\Sol$ is just
the standard volume on $\Ra^3$. In the case when $m \neq n$ we
introduce a new measure. In the case of $\Sol$ this is just
$\mu=dxdydz$.  Note that on $z$ level sets this is a rescaling of
$\vol$ by a factor of $e^{(n-m)z}$. Analogously on $\DL(m,n)$, we
choose a height function $h:\DL(m,n){\rightarrow}\Za$ and let
$\mu$ be counting measure times $n^{h(x)}m^{-h(x)}$. Recall that
we are assuming that $m{\geq}n$.  The measure $\mu$ is also
introduced in \cite{BLPS} and is natural for many problems.

\medskip
We now define certain useful subsets of $\Sol$. We define these sets
simply as subsets of $\Ra^3$. Let $ B(L,\vec{0})
= [-\frac{e^{2mL}}{2},\frac{e^{2mL}}{2}] \times
[-\frac{e^{2nL}}{2},\frac{e^{2nL}}{2}] \times
[-\frac{L}{2},\frac{L}{2}]$. When $m=n$, then $|B(L,\vec{0})|
\approx L e^{2mL}$ and $Area(\partial B(L,\vec{0})) \approx
e^{2mL}$, so $B(L)$ is a F\"olner set.

To define the analogous object in $\DL(m,n)$, we look at the set
of points in $\DL(m,n)$ we fix a basepoint $(\vec{0})$ and a
height function $h$ with $h(\vec{0})=0$. Let $L$ be an even
integer and let $\DL(m,n)_L$ be the
$h{\inv}([-\frac{L+1}{2},\frac{L+1}{2}])$. Then $B(L,\vec{0)})$ is
the connected component of $\vec{0}$ in $\DL(m,n)_L$.  We are
assuming that the top and bottom of the box are midpoints of
edges, to guarantee that they have zero measure.

 We call $B(L,\vec{0})$ a box of size $L$
centered at the identity. In $\Sol$, we define the box of size $L$
centered at a point $p$ by $B(L,p)=T_pB(L,\vec{0})$ where $T_p$ is
left translation by $p$. We frequently omit the center of a box in
our notation and write $B(L)$.  For the case of $\DL(m,n)$ it is
easiest to define the box $B(L,p)$ directly. That is let
$\DL(m,n)_{[h(p)-\frac{L+1}{2},h(p)+\frac{L+1}{2}]}=h^{\inv}([h(p)-\frac{L+1}{2},h(p)+\frac{L+1}{2}])$
and let $B(L,p)$ be the connected component of $p$ in
$\DL(m,n)_{[h(p)-{\frac{L+1}{2}},h(p)+\frac{L+1}{2}]}$.  It is
easy to see that isometries of $\DL(m,n)$ carry boxes to boxes.

We record the following lemma which holds for any model space
$\X$.

\begin{lemma}
\label{lemma:box:amenable} When $m=n$, the fraction of the volume
of $B(L)$ which is within $\epsilon L$ of the boundary of $B(L)$
is $O(\epsilon)$. In all other cases, this is true for the
$\mu$-measure but not the volume.
\end{lemma}

\mc{ Do we prove this? It is rather clear.}

We first describe $B(L)$ in the case of $\Sol$. In this case, the
top of $B(L)$, meaning the set
$[-\frac{e^{2mL}}{2},\frac{e^{2mL}}{2}] \times
[-\frac{e^{2nL}}{2},\frac{e^{2nL}}{2}] \times \{\frac{L}{2}\}$, is
not at all square - the sides of this rectangle are horocyclic
segments of lengths $e^{2mL}$ and $1$ - in other words it is just
a small metric neighborhood of a horocycle. Similarly, the bottom
is also essentially a horocycle but in the transverse direction.
Further, we can connect the $1$-neighborhood of  any point of the
top horocycle to the $1$-neighborhood of any point of the bottom
horocycle by a vertical geodesic segment, and these segments
essentially sweep out the box  $B(L)$.   This picture is even
easier to understand in the Diestel-Leader graphs $\DL(n,n)$,
where the boundary of the box is simply the union of the top and
bottom "horocycles", and the vertical geodesics in the box form a
complete bipartite graph between the two. Thus a box $B(L)$
contains a very large number of quadrilaterals.

\subsection{\bf Discretizing Sol}

We describe in this section a variety of ways of seeing more closely
the analogy between the geometry of $\Sol$ and the geometry of
$\DL(m,n)$.  This is done most easily by thinking about
discretizations of $\Sol$.  While we do not use these
discretizations formally in our proof, they are the reason why we
sometimes only describe a proof completely in one of the model
geometries.

To see this picture most clearly, we first remark that in a box
$B(R)$ in $\DL(m,n)$, one can form an auxiliary graph $\hat B(R)$
whose vertices consist only of those vertices on the top and bottom
of the graph and where there is an edge between vertices whenever
there is a vertical geodesic connecting them.  This graph is
complete bipartite, where the parts are the top and bottom of the
box.

In $\Sol$ one can make a similar construction.  Namely given $B(R)$,
we construct a graph $\hat B(R)$ as follows.  Choose a $C$-net in
the top and bottom of the box.  Vertices will be the points in the
$C$-net with the bipartition into those on the top and those on the
bottom.  Connect a vertex $x$ on the top to a vertex $y$ on the
bottom if there is a vertical leaving the $10C$ neighorhood of $x$
arriving in the $10C$ neighborhood of $y$. (The constants $C$ and
$10$ are arbitrary.)  It is an elementary exercise in hyperbolic
geometry to show that $\hat B(R)$ is complete bipartite.

While we do not use the graphs $\hat B(R)$ explicitly in this paper,
they contain much of the geometry that is necessary for our
arguments.

\subsection{\bf Tiling}
\label{sec:tiling}

The purpose of this subsection is to prove the following lemma.

\begin{lemma}
\label{lemma:tiling} Choose constants $L>R$ such that $L/R{\in}\Za$. We can write
\begin{equation}
\label{eqn:tiling} B(L)=\bigsqcup_{i \in I} B_i(R) \sqcup \Upsilon
\end{equation}
 where
$\mu(\Upsilon)=O(R/L)\mu(B(L))$ and the implied constant depends
only on the model space.

In the case when $\X=\DL(m,n)$, then $\Upsilon$ can be chosen to
be empty.  (This is also possible for $\Sol$ if $e^m$ and $e^n$
are integers.)
\end{lemma}

{\noindent}{\bf Remark:} We will always refer to a decomposition
as in equation (\ref{eqn:tiling}) as a {\em tiling} of $B(L)$. We
often omit specific reference to the set $\Upsilon$ when
discussing tilings.

\bold{Proof.} For simplicity of notation, we assume $B(L)$ is
centered at the origin.

We give the proof first in the case of $\DL(m,n)$ where it is
almost trivial. Since $\frac{L}{R}{\in}\Za$, we can partition
$[-\frac{L+1}{2},\frac{L+1}{2}]$ into subsegments of length $R$
which we label $S_1, \cdots, S_{J}$ where $J=\frac{L}{R}$.  We can
then look at the sets $\DL(m,n)_j=h{\inv}(S_j)$. Each connected
component of $\DL(m,n)_j$ is clearly a box $B_{j,k}(R)$ of size
$R$. Each $B_{j,k}(R)$ is either entirely inside or entirely
outside of $B(L)$.  We choose only those $k$ for which $B_{j,k}(R)
\subset B(L)$. It is also clear that, after reindexing, we have
chosen boxes such that $B(L)= \bigsqcup B_i(R)$.

In $\Sol$ the proof is similar, though does not in general give an
exact tiling.  We simply take the box $B(L)$ and cover it as best
possible with boxes of size $R$.  Since $\frac{L}{R}{\in}\Za$, if
we take $B(R,\vec{0})$ and look at translates by $(0,0,Rc)$ for
$c$ an integer between $-\frac{L}{R}$ and $\frac{L}{R}$, the
resulting boxes are all in $B(L)$.  We then take the resulting box
$B(R, (0,0,Rc))$ at height $k$ and translate it by vectors of the
form $(\vec{ae^{mcR},be^{ncR},0})$ where $  |a| \leq
e^{m(L-(c+1)R)} $ and $|b| \leq e^{n(L-(1-c)R)}$ are integers.
This results in boxes $B(R,(a,b,c))$ which we re-index as
$B_i(R)$.  It is clear that every point not in $\bigsqcup _i
B_i(R)$ is within $R$ of the boundary of $B(L)$.  Letting
$\Upsilon=B(L) - \bigsqcup _i B_i(R)$ we have that
$\mu(\Upsilon)<O(R/L)\mu(B(L))$ by Lemma \ref{lemma:box:amenable}.

\qed\medskip

\section{Step 1}

All the results of this section hold for $\X$ with $m{\geq}n$, so
in particular for the case $m=n$.  The case $m=n$ will be used in
the sequel \cite{EFW1}.  Also, all results in this section hold
for quasi-isometric embeddings, i.e. maps satisfying $(1)$ but not
$(2)$ of Definition \ref{defn:qi}.  Before stating the main
result of this part of the paper, we make some definitions.  The first
is simple and just says that a map is close to a product map, where here
close depends on the diameter of the domain of definition.

\begin{defn}
\label{defn:epsilonsublineartostandard}
Let $E$ be subset of $\Sol$ of diameter $R$. A quasi-isometric embedding 
$\phi: E \rightarrow \Sol$ is called {\em $\epsilon$-sublinear to a product map}
if there is product map $\hat \phi$ from $\Sol$ to $\Sol$ such that 
$d(\hat \phi|_E, \phi) \leq O(\epsilon R)$
\end{defn}

Our arguments would be much simpler if we could show quickly that $\phi$ restricted to a box $B(R)$ was
$\epsilon$-sublinear to a product.  The weaker statement that we prove in this section requires another
definition.

\begin{defn}
\label{defn:thetamostly}
Given constants $R<L$, a box $B(L)$ and a quasi-isometry $\phi:B(L) \rightarrow \Sol$, we say that $\phi$ is
{\em $\theta$ mostly $\epsilon$-sublinear to product maps at scale $R$} if one can tile
\begin{displaymath}
B(L) = \bigsqcup_{i \in I} B_i(R)
\end{displaymath} 
and there exists a subset $I_g$ of $I$ with $\mu(\bigcup_{i \in
I_g} B_i(R))  \ge (1-\theta) \mu(B(L))$ so that for any $i \in
I_g$ there exists $U_i \subset B_i(R)$ with $\mu(U_i) \ge
(1-\theta) \mu(B_i(R))$ such that $\phi$ restricted to each $U_i$
is $\epsilon$-sublinear to a product map.
\end{defn}

When considering maps that are $\theta$ mostly $\epsilon$-sublinear to product maps at scale $R$, 
we will denote by $\hat \phi_i$ the product map that is $\epsilon$-sublinear to $\phi$ on $U_i$.  Note that
the definition allows $\hat \phi_i \neq \hat \phi_j$.

In the this part of the paper,
our aim is to prove the following:

\begin{theorem}
\label{theorem:main:step:one} Suppose $\theta > 0$, $\epsilon >
0$. Then there exist constants $0 < \alpha < \beta < \Delta$
(depending on $\theta$, $\epsilon$, $\kappa$, $C$ and the model
spaces) such that the following holds: Let $\phi: \X \to \pX$ be a
$(\kappa,C)$ quasi-isometry and suppose $r_0$ is sufficiently
large (depending on $\kappa$, $C$, $\theta$, $\epsilon$). Then for
any $L > \Delta r_0$ and any $B(L)$, there exists $R$ with $\alpha
r_0 < R < \beta r_0$ such that $\phi$ is $\theta$ mostly $\epsilon$-sublinear to
product maps scale $R$.
\end{theorem}

\bold{Remarks.} This theorem says that every sufficiently large
box $B(L)$ can be tiled by much smaller boxes $B_i(R)$, and for
most (i.e. $1-\theta$ fraction) of the smaller boxes $B_i(R)$
there exist a subset $U_i$ containing $1-\theta$ fraction of the
$\mu$-measure of $B_i(R)$ on which the map is a product map, up to
error $O( \epsilon R) \ll R$.  In the case where $n=m$, the
measure of $B_i(R)$ is independent of $i$, and we have exactly
$|I_g| \geq (1-\theta)|I|$.  When $m{\neq}n$, both the number of
boxes of size $R$ in a height level set tiling and $\mu(B_i(R))$
are functions of height.  We note here that it is possible to
apply the proof of Theorem \ref{theorem:main:step:one}
simultaneously to a finite collection $J$ of boxes $B_j(L)$ all of
the same size and obtain the same conclusions (with the same
constants) on most of the boxes in $J$.  As long as $m=n$, by most
boxes in $J$ we mean most boxes with the counting measure on $J$.
This observation will be used in \cite{EFW1}.

One should note that the number $R$, and the subset where we
control the map, depends on $\phi$. Also in
Theorem~\ref{theorem:main:step:one} there is no assertion that the
product maps $\hat{\phi}_i$ on the different boxes $B_i(R)$ match
up.

This theorem is in a sense an analogue of Rademacher's theorem
that a lipschitz function (or map) is differentiable almost
everywhere. The boxes $B_i(R)$ with $i{\in}I_g$ should be thought
of as coarse analogues of points of differentiability.

 The proof of this theorem is done in several steps. First
we apply a {\em coarse differentiation} argument to show that there
exist $R,L,I_g$ etc. as in Theorem~\ref{theorem:main:step:one} such
that restricted to each $B_i(R)$ with $i{\in}I_g$, the map $\phi$
sends most vertical geodesics to within $O(\epsilon R)$ of a
vertical geodesic.  In the second step, we use some elementary
geometry of the model space and particularly of the set $B_i(R)$ to
show that this implies that $\phi$ is close to a product map on most
of the measure of $B_i(R)$.  In particular, we apply Lemma
\ref{lemma:quadrilateral} in the range of $\phi$ to the images of
quadrilaterals in the domain to ensure these images have essentially
the same geometric structure.  That this is enough to control the
map on $B_i(R)$ essentially follows from the fact that $B_i(R)$ is
basically a complete bipartite graph on the top and bottom of the
box.

\subsection{Behavior of quasi-geodesics}
\label{sec:quasigeodesics}

We begin by discussing some quantative estimates on the behavior
of quasi-geodesic segments in $\pX$ (or equivalently in $\X$).
Throughout the discussion we assume $\alpha:[0,r] \to \pX$ is a
$(\kappa,C)$-quasi-geodesic segment for a fixed choice of
$(\kappa,C)$, i.e. $\alpha$ is a quasi-isometric embedding of
$[0,r]$ into $\pX$.  A quasi-isometric embedding is a map that
satisfies point $(1)$ in Definition \ref{defn:qi} but not point
$(2)$. All of our quasi-isometric embeddings are assumed to be
continuous.

\begin{defn}[{\bf $\epsilon$-monotone}]
A quasigeodesic segment $\alpha: [0,r] \to
  \pX$ is  $\epsilon$-monotone  if for all
  $t_1, t_2 \in [0,r]$ with $h(\alpha(t_1)) =
  h(\alpha(t_2))$ we have $|t_1 - t_2 | < \epsilon r$.
\end{defn}

\begin{figure}[ht]
\begin{center}
\includegraphics[width=1in]{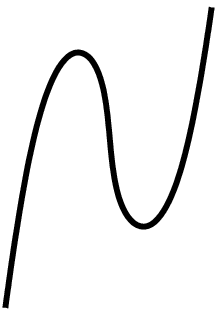}
\caption{A quasigeodesic segment which is not
$\epsilon$-monotone.} \label{fig:not_monotone}
\end{center}
\end{figure}

In \S\ref{sec:staircase} and \S\ref{sec:fullmeasure} we will also
need a variant.  The reader may safely ignore this variant on
first reading this section.

\begin{definition}[Weakly $(\eta,C_1)$-monotone]
\label{def:weakly:eta:monotone} A quasigeodesic segment $\alpha:
[0,r] \to
  \pX$ is {\em weakly $(\eta,C_1)$-monotone} if for any two points $0 < t_1
  < t_2 < r$ with $h(\alpha(t_1)) = h(\alpha(t_2))$, we
  have $t_2 - t_1 < \eta t_2+C_1$.
\end{definition}

{\noindent}{\bf Remark}: An $\epsilon$-monotone quasi-geodesic
$\alpha:[0,r]{\rightarrow}\pX$ is a weakly $(\epsilon,\epsilon
r)$-monotone quasi-geodesic.

The following fact about $\varepsilon$-monotone geodesics is an
easy exercise in hyperbolic geometry:

\begin{lemma}
\label{lemma:epsilon:monotone:segment} $ $
\begin{itemize}
\item[{\rm (a)}] Suppose $\alpha: [0,r] \to
  \pX$ is an $\epsilon$-monotone quasi-geodesic segment. Then, there
  exists a vertical geodesic segment $\lambda$ in $\pX$ such that
  $d(\alpha, \lambda) \le \omega_1 \epsilon r$, where $\omega_1$
  depends only on the model space $\pX$.
\item[{\rm (b)}] Suppose $\alpha: [0,r] \to
  \pX$ is a weakly  $(\eta,C_1)$-monotone quasi-geodesic segment.
  Then, there
  exists a vertical geodesic segment $\lambda$ in $\pX$ such that
  $d(\bar{\gamma}(t), \lambda(t)) \le 2 \kappa \eta t + \omega_2 C_1$,
  where $\omega_2$ depends only on $\pX$.
\end{itemize}
\end{lemma}

\bold{Proof.} Both $\epsilon$-monotone and weakly
$(\eta,C)$-monotone imply that the projections of $\alpha$ onto
both $xz$ and $yz$ hyperbolic planes are quasi-geodesics.
 The result is then a consequence of the Mostow-Morse lemma and the
 fact that the only geodesics shared by both families of hyperbolic
 planes are vertical geodesics.  One can also prove the lemma by
 direct computation.
\qed\medskip

\noindent{\bf Remark:} The distance $d(\alpha,\lambda)$ in (a) is
the Hausdorff distance between the sets and does not depend on
parametrizations. However, the parametrization on $\lambda$
implied in (b) is not neccessarily by arc length. \medskip

\begin{lemma}[Subdivision]
\label{lemma:subdivision} Suppose $\alpha: [0,r] \to \pX$ is a
quasi-geodesic segment which is not $\epsilon$-monotone and $r \gg
C$. Suppose $N \GG 1$ (depending on $\epsilon$, $\kappa$, $C$).
Then
\begin{displaymath}
\sum_{j=0}^{N-1}\left| h(\alpha(\tfrac{(j+1)r}{N})) -
  h(\alpha(\tfrac{jr}{N})) \right| \ge
\left|h(\alpha(0)) - h(\alpha(r))\right| + \frac{\epsilon r}{ 8
\kappa^2}.
\end{displaymath}
\end{lemma}

\begin{figure}[ht]
\begin{center}
\includegraphics[width=1.0in]{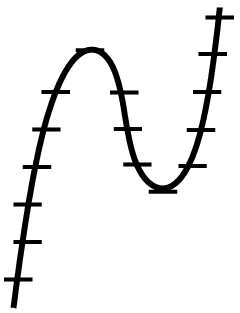}
\caption{Proof of Lemma~\ref{lemma:subdivision}}
\label{fig:subdivision}
\end{center}
\end{figure}

Informally, the proof amounts to the assertion that if $N$ is
sufficiently large, the total variation of the height increases
after the subdivision by a term proportional to $\epsilon$. See
Figure~\ref{fig:subdivision}.
\medskip

\bold{Proof.} Without loss of generality, we may assume that
$h(\alpha(0)) \ge h(\alpha(t_1)) = h(\alpha(t_3)) \ge
h(\alpha(r))$, where $0 = t_0 < t_1 < t_3 < t_4 = r$ (if not,
parametrize in the opposite direction). Since $t_3 -t_1 > \epsilon
r$, $\alpha(t_3)$ and $\alpha(t_1)$ are two points in $\pX$ which
are at the same height and are at least $\epsilon r/\kappa$ apart.
Then, by $\pX$ geometry, any long enough $(\kappa,C)$-quasigeodesic path connecting
$\alpha(t_3)$ and $\alpha(t_1)$ must contain a point $q$ such that
$|h(q) -  h(\alpha(t_1))| \ge (\epsilon r)/(4\kappa)$. Hence,
there exists a point $t_2$ with $t_1 < t_2 < t_3$ such that
$|h(\alpha(t_2)) -  h(\alpha(t_1))| \ge (\epsilon r)/(4\kappa)$.
Hence,
\begin{displaymath}
\sum_{j=1}^4 |h(\alpha(t_j)) - h(\alpha(t_{j-1}))| \ge
|h(\alpha(0)) - h(\alpha(r))| + \frac{\epsilon r}{4 \kappa}
\end{displaymath}
If $N$ is large enough then the points $t_1$, $t_2$ and $t_3$ have
good approximations of the form $j r/N$, with $j \in \zed$. This
implies the lemma. \qed\medskip

\noindent {\bf Choosing Scales:} Choose $1 \ll r_0 \ll r_1 \ll
\dots \ll r_S$. In particular, $C \ll r_0$ and for $s \in [0,S-1]
\cap \zed$, $r_{s+1}/r_s > N$ where $N$ is as in Lemma
\ref{lemma:subdivision}.

\begin{lemma}
\label{lemma:scales} Suppose $L \GG r_S$, and suppose $\alpha: [0,
L] \to \pX$ is a quasi-geodesic segment. For each $s \in [1,S]$,
subdivide $[0,L]$ into $L/r_s$ segments of length $r_s$. Let
$\delta_s(\alpha)$ denote the  fraction  of these segments whose
images are not $\epsilon$-monotone. Then,
\begin{displaymath}
\sum_{s=1}^S \delta_s(\alpha) \le \frac{16\kappa^3}{\epsilon}.
\end{displaymath}
\end{lemma}

\noindent{\bf Remark:} The utility of the lemma is that the right
hand side is fixed and does not depend on $S$.  So for $S$ large
enough, some (in fact many) $\delta_s(\alpha)$ must be small.

\noindent {\it Proof.} By applying Lemma~\ref{lemma:subdivision}
to each non-$\epsilon$-monotone segment on the scale $r_S$, we get
\begin{multline*}
\sum_{j=1}^{L/r_{S-1}} \left|h(\alpha(j r_{S-1})) -
  h(\alpha((j-1)r_{S-1})) \right| \ge \\ \ge
\sum_{j=1}^{L/r_S} \left|
  h(\alpha(j r_S)) -
  h(\alpha((j-1)r_S))\right|+\delta_S(\alpha)
\frac{\epsilon L}{8 \kappa^2}.
\end{multline*}
Doing this again, we get after $S$ iterations,
\begin{multline*}
\sum_{j=1}^{L/r_0} \left|h(\alpha(j r_0)) -
  h(\alpha((j-1)r_0)) \right| \ge \\ \ge
\sum_{j=1}^{L/r_S} \left|
  h(\alpha(j r_S)) -
  h(\alpha((j-1)r_S))\right|+
\frac{\epsilon L}{8 \kappa^2} \sum_{s=1}^S \delta_s(\alpha).
\end{multline*}
But the left-hand-side is bounded from above by the length and so
bounded above by $2\kappa L$. \qed\medskip

\subsection{Averaging}
\label{sec:averaging}

In this subsection we apply the estimates from above to images of
geodesics under a quasi-isometry from $\X$ to $\pX$. The idea is
to average the previous estimates over families of geodesics.  In
order to unify notation for the two possible model space types, we
shift the parametrization of vertical geodesics in $\DL(m,n)$ so
that they are parametrized by height minus $\frac{1}{2}$, i.e. by
the interval $[-\frac{L}{2},\frac{L}{2}]$ rather than
$[-\frac{L+1}{2},\frac{L+1}{2}]$.

\bold{Setup and Notation.}
\begin{itemize}
\item Suppose $\phi: \X \to \pX$ is a $(\kappa,C)$ quasi-isometry.
Without loss of generality, we may assume that $\phi$ is
continuous. \item Let $\gamma: [-\frac{L}{2},\frac{L}{2}] \to \X$
be a vertical geodesic segment parametrized by arclength where $L
\gg C$. \item Let $\overline{\gamma} = \phi \circ \gamma$. Then
$\overline{\gamma}: [-\frac{L}{2},\frac{L}{2}] \to \pX$ is a
quasi-geodesic segment.
\end{itemize}

\noindent It follows from Lemma~\ref{lemma:scales}, that for every
$\theta > 0$ and every geodesic segment $\gamma$, assuming that
$S$ is sufficiently large, there exists $s \in [1,S]$ such that
$\delta_s(\overline{\gamma}) < \theta$. The difficulty is that $s$
may depend on $\gamma$. In our situation, this is overcome as
follows:
\medskip


We will average the result of Lemma~\ref{lemma:scales} over $Y_L$,
the set of vertical geodesics in $B(L)$. Let $|Y_L|$ denote the
measure/cardinality of $Y_L$. We will always denote our average by
$\Sigma$, despite the fact that when $\X=\Sol$ this is actually an
integral over $Y_L$ and not a sum.  When $\X=\DL(m,n)$ it is
actually a sum.  Changing order, we get:
\begin{displaymath}
\sum_{s=1}^S \left( \frac{1}{|Y_L|} \sum_{\gamma \in Y_L}
    \delta_s(\overline{\gamma}) \right) \le \frac{16 \kappa^3}{\epsilon}.
\end{displaymath}
Let $\delta > 0$ be a small parameter (In fact, we will choose
$\delta$ so that $\delta^{1/4} = \min(\epsilon,\theta/256)$, where
$\theta$ is as in Theorem~\ref{theorem:main:step:one}).  \mc{
check} Then, if we choose $S > \frac{16 \kappa^3}{\epsilon
\delta^4}$, then there exists a scale $s$ such that
\begin{equation}
\label{eq:average:delta:small} \frac{1}{|Y|}\sum_{\gamma \in Y}
  \delta_s(\bar{\gamma})  \le \delta^4.
\end{equation}

\bold{Conclusion.} On the scale $R \equiv {r_s}$, at least
$1-\delta^4$ fraction of all vertical geodesic segments of length
$R$ in $B(L)$
have nearly vertical images under $\phi$. 
\medskip

\relax From now on, we fix this scale, and drop the index $s$. We
will refer to segments of length $R$ arising in our subdivision as
{\em edges} of length $R$. In the case of $DL(m,n)$ these edges
are unions of edges in the graph.  In what follows we will use the
terms {\em big edges} for edges of length $R$ if there is any
chance of confusion with an actual edge in the graph $DL(m,n)$.

\bold{Remark.} The difficulty is that, at this point, even though
we know that most edges have images under $\phi$ which are nearly
vertical, it is possible that some may have images which are going
up,  and some may have images which are going down.

\subsection{Alignment}
\label{sec:align}

We assume that $L/R \in \zed$. As described in \S\ref{sec:tiling},
we tile
\begin{displaymath}
B(L) = \bigsqcup_{i \in I} B_i(R).
\end{displaymath}
Let $Y_i$ denote the set of vertical geodesic segements in
$B_i(R)$. We have
\begin{equation}
\label{eq:tiling:identity} \frac{1}{|Y_L|}\sum_{\gamma \in Y_L}
  \delta_s(\bar{\gamma})  = \sum_{i \in I} \frac{\mu(B_i(R))}{\mu(B(L )}\left( \frac{1}{|Y_i|}
  \sum_{\lambda \in Y_i} \delta_s({\bar{\lambda}})\right) + O(\tfrac{R}{L}),
\end{equation}
where $\bar{\lambda} = \phi \circ \lambda$, and
$\delta(\bar{\lambda}) = \delta_s(\bar{\lambda})$ is equal to $0$
if $\lambda$ is $\epsilon$-monotone, and equal to $1$ otherwise.
The error term of $O(R/L)$ is due to the fact that the tiling may
not be exact, see Lemma~\ref{lemma:tiling}. To justify equation
\ref{eq:tiling:identity}, one uses that $\mu(B(L))=|Y_L|L$ and
$\mu(B_i(R))=|Y_i|R$.

Since the left hand side is bounded by $\delta^4$ and assuming $R/L
\ll \delta^2$, we conclude the following:


\begin{lemma}

\label{lemma:good:box} Let us tile $B(L)$ by boxes $B_i(R)$ of
size $R$, so that $B(L) = \bigsqcup_{i \in I} B_i(R)$. Then there
exists a subset $I_g$ of the indexing set $I$ with $\mu(\bigcup_{i
\in I_g} B_i(R))  \ge (1-\delta^2) \mu(B(L))$ such that if we let
$Y_i$ denote the set of vertical geodesics in $B_i(R)$ then
\begin{equation}
\label{eq:goodbox:average:delta:small} \frac{1}{|Y_i|}\sum_{\gamma
\in Y_i}
  \delta_s(\bar{\gamma}) \le 2 \delta^2.
\end{equation}
\end{lemma}

\noindent Note that $R$ is the length of one big edge so that the
set $Y_i$ of vertical geodesics in $B_i(R)$ consists of big edges
connecting the top to the bottom. The equation
(\ref{eq:goodbox:average:delta:small}) means that the fraction
these edges which are not $\epsilon$-monotone is at most $2
\delta^2$.

\bold{Notation.} In the rest of \S\ref{sec:align} and in
\S\ref{sec:product}  we fix $i \in I_g$ and drop the index $i$.
We refer to a vertical geodesic segment $e$ running from bottom
to top of $B(R)$ as an {\em edge} of $B(R)$.  We say that
$e$ is ``upside-down'' if $\phi(e)$ is going down, and
``right-side-up'' if $\phi(e)$ is going up.

\medskip\noindent

\begin{lemma}[Alignment]
\label{lemma:big:edges:aligned} Let $e$ be an $\epsilon$-monotone
big edge of $B(R)$ going from the bottom to the top. Then either the
fraction of the big edges in $B(R)$ which are upside down or the
fraction of the big edges in $B(R)$ which are right-side-up is at
least $1-4 \delta$.
\end{lemma}

\bold{Proof.} We have a natural notion of ``top'' vertices and
``bottom'' vertices so that each big edge connects a bottom vertex
to a top vertex.  Then $B(R)$ is a complete bipartite graph. There
must be a subset $E$ of vertices of density $1-4\delta$, such that
for each vertex in $v \in E$ the fraction of the edges incident to
$v$ which are not $\epsilon$-monotone is at most $\delta/2$. Let
$\Gamma_1$ be the subgraph of $B(R)$ obtained by erasing any edge
$e$ such that $\phi(e)$ is not $\epsilon$-monotone. We orient each
edge $e$ of $\Gamma_1$ by requiring that $\phi(e)$ is going down.

Let $p_1, p_2 \in E$ be any two top vertices in the good set. Then we can find two
bottom vertices $q_1, q_2$ such that all four quasigeodesic
segments $\phi(\overline{p_1 q_1})$, $\phi(\overline{p_1 q_2})$,
$\phi(\overline{p_2 q_1})$ and $\phi(\overline{p_2 q_2})$ are all
$\epsilon$-monotone, $\overline{p_1 q_1}$ and $\overline{p_1 q_2}$
diverge quickly at $p_1$, and $\overline{p_2 q_1}$ and
$\overline{p_2
  q_2}$ diverge quickly at $p_2$.  We can arrange for the fast divergence, since
  fast divergence occurs generically, i.e. on the complement of a set of small measure.


We now apply Lemma~\ref{lemma:quadrilateral} to conclude that
$h(\phi(p_1)) = h(\phi(p_2))+ O(\epsilon R)$ and that all the
segments $\phi(\overline{p_i q_j})$ with $i,j=1,2$ have the same
orientation. Thus any two top vertices in $E$ have images on
essentially the same height, say $h_1$. Similarly, any two bottom
vertices in $E$ have images on the same height, say $h_2$. Since
we must have $h_1
> h_2$ or $h_1 < h_2$, the lemma holds. \qed\medskip

We define the {\em dominant orientation} to be right-side-up or
upside-down so that the fraction of big edges which have the
dominant orientation is at least $1 - 4 \delta$.

\subsection{Construction of a product map}
\label{sec:product}

 Recall that $Y$ is the set of vertical geodesics in $B(R)$.
Let $Y'$ denote the space of pairs $(\gamma,x)$ where $\gamma \in
Y$ is a vertical geodesic in $B(R)$ and $x \in \gamma$ is a point.
Let $| \cdot |$ denote uniform measure on $Y'$. (In the case of
$\DL(m,n)$ this is just the counting measure.)  The following
lemma is a formal statement regarding subsets of $Y'$ of large
measure.

\begin{lemma}
\label{lemma:uniform:set} Suppose $R \GG 1/\theta_1$ (where the
implied constant depends only on the model space). Suppose $E
\subset Y'$, with $|E| \ge (1 - \theta_1) |Y'|$. Then, there
exists a subset $U \subset B(R)$ such that:
\begin{itemize}
\item[{\rm (i)}] $\mu(U) \ge (1-2\sqrt{\theta_1}) \mu(B(R))$,
where
  $\mu$ is defined in $\S\ref{subsection:volumeandmeasure}$.
\item[{\rm (ii)}] If $x \in U$, then for at least
  $(1-\sqrt{\theta_1})$ fraction of the vertical
  geodesics $\gamma \in Y$ passing within distance $1/2$ of $x$,
  $(\gamma,x) \in E$.
\end{itemize}
\end{lemma}

\bold{Remark.} Note that for the case of $\DL(m,n)$, any geodesic
passing within distance $(1/2)$ of $x$ passes through $x$.

\bold{Proof.} For $x \in B(R)$, let $Y(x) \subset Y$ denote the
set of geodesics which pass within $1/2$ of $x$. For clarity, we
first give the proof for the $\DL(m,n)$ case.  Note that $|Y(x)| =
c \mu(\{x\})$, where $c$ depends on $m,n$ and the location and
size of $B(R)$. Note that
\begin{equation}
\label{eq:measure:exchange} |Y| R = \sum_{x \in B(R)} \sum_{\gamma
\in Y(x)}  1 = \sum_{x \in
  B(R)} |Y(x)| = \sum_{x \in B(R)} c \mu(\{x\}) = c \mu(B(R))
\end{equation}

Suppose $f(\gamma,x)$ is any function of a geodesic $\gamma$ and a
point $x \in \gamma$: then,
\begin{align}
\label{eq:exchange} \frac{1}{|Y|R} \sum_{\gamma \in Y} \sum_{x \in
\gamma} f(\gamma, x) & = \frac{1}{|Y| R} \sum_{x \in B(R)}
\sum_{\gamma \in Y(x)} f(\gamma,x) \notag \\
& = \frac{1}{|Y| R} \sum_{x \in B(R)} \frac{1}{|Y(x)|}
\sum_{\gamma \in Y(x)} |Y(x)| f(\gamma,x) \notag \\ & =
\frac{1}{\mu(B(R))} \sum_{x \in B(R)} \frac{1}{|Y(x)|}
\sum_{\gamma \in Y(x)} \mu(\{x\}) f(\gamma,x),
\end{align}
where in the last line we used (\ref{eq:measure:exchange}).

We apply (\ref{eq:exchange}) with $f$ the the characteristic
function of the complement of $E$. We get,
\begin{equation}
\label{eq:chi} \frac{1}{\mu(B(R))} \sum_{x \in B(R)} \mu(\{x\})
\left(\frac{1}{|Y(x)|} \sum_{\gamma \in Y(x)}  f(\gamma,x) \right)
< \theta_1
\end{equation}
Let $F(x)$ denote the parenthesized quantity in the above
expression. Let $E_2 = \{x \in B(R) \st F(x) > \sqrt{\theta_1}
\}$. Recall that Markov's inequality says that for any real-valued
function $f$, and any real number $a > 0$, the measure of the set
$\{ |f| > a \}$ is at most $\frac{1}{a} \int |f|$. Then, by this
inequality, $\mu(E_2)/\mu(B(R)) \le \sqrt{\theta_1}/\theta_1 =
\sqrt{\theta_1}$, and for $x \not\in E_2$, for at least
$(1-\sqrt{\theta_1})$ fraction of the geodesics $\gamma$ passing
through $x$, $(\gamma,x) \in E$.

This completes the proof for the $\DL(m,n)$ case. In $\Sol$ the
computation is essentially the same, except for the fact that
$|Y(x)|$ (i.e. the measure of set of geodesics passing within
$(1/2)$ of $x$) can become smaller when $x$ is within $(1/2)$ of
the boundary of $B(R)$. However, the relative $\mu$ measure of
such points is $O(1/R)$ by Lemma~\ref{lemma:box:amenable}.
Therefore (\ref{eq:measure:exchange}) and (\ref{eq:exchange}) hold
up to error $O(1/R) < \theta_1$. \qed\medskip


\begin{corollary}
\label{cor:mostly:monotone} There exists a subset $U \subset B(R)$
with $\mu(U) > (1 - 8\sqrt{\delta}) \mu(B(R))$ such that for $x
\in U$, $(1-2 \sqrt{ \delta})$-fraction of the geodesics passing
within $(1/2)$ of $x$ have  $\epsilon$-monotone image under $\phi$  and have images
with the dominant orientation.
\end{corollary}

\bold{Proof.} Let $E$ denote the set of pairs $(\gamma,x)$ where
$\gamma \in Y$ is a dominantly oriented $\epsilon$-monotone
geodesic segment, and $x$ is a point of $\gamma$. Let $U \subset
B(R)$ be the subset constructed by Lemma~\ref{lemma:uniform:set}.
Since $|E| \ge (1 - 4\delta) |Y'|$, $\mu(U) \ge (1 -8
\sqrt{\delta}) \mu(B(R))$. \qed\medskip

\begin{lemma}
\label{lemma:height:preserving} Suppose $\phi$ and $B(R)$ and $U$
are as in Corollary~\ref{cor:mostly:monotone}.  Then, there exist
functions $\psi:\Ra^3{\rightarrow}\Ra^2$, $q:\Ra{\rightarrow}\Ra$,
and a subset $U_1 \subset B(R)$ with $\mu(U_1) > (1 - 128
\delta^{1/4}) \mu(B(R))$ \mc{ CHECK} such that for $(x,y,z) \in
U_1$,
\begin{equation}
\label{eq:partial:product} d(\phi(x,y,z),(\psi(x,y,z),q(z))) =
O(\epsilon R)
\end{equation}
\end{lemma}

\bold{Proof.} We assume that the dominant orientation is
right-side-up (the other case is identical). Now suppose $p_1, p_2
\in U$ belong to the same $x$-horocycle. By the construction of
$U$ there exist $q_1, q_2$ in $B(R)$ (above $p_1$, $p_2$) such
that for each $i=1,2$ the two geodesic segments $\overline{p_i
  q_1}$ and $\overline{p_i q_2}$ leaving $p_i$ diverge quickly, and
 each of the
quasigeodesic segments $\phi(\overline{p_i q_j})$ is
$\epsilon$-monotone. Then by
Lemma~\ref{lemma:epsilon:monotone:segment}, each of the
$\phi(\overline{p_i q_j})$ is within $O(\epsilon R)$ of a
quasi-geodesic segment $\lambda_{ij}$. Now by applying
Lemma~\ref{lemma:quadrilateral} to the $\lambda_{ij}$ we see that
$\phi(p_1)$ and $\phi(p_2)$ are on the same $x$-horocycle, up to
an error of $O(\epsilon R)$. Thus, the restriction of $\phi$ to
$U$ preserves the $x$-horocycles. A similar argument (but now we
will pick $q_1$, $q_2$ below $p_1$, $p_2$) shows that the
restriction of $\phi$ to $U$ preserves the $y$-horocycles. We can
now conclude that $\phi$ is height respecting
on a slightly smaller set $U_1$ \mc{ NEW LEMMA} \mc{
need to justify $\psi$ not depending on $z$}, i.e. there exist
functions $\psi: \reals^3 \to \reals^2 $ and $q: \reals \to
\reals$ such that for $(x,y,z) \in U_1$,
(\ref{eq:partial:product}) holds. \qed\medskip

\begin{proposition}
\label{prop:product:map} Suppose $\phi$ and $B(R)$ and $U$ are as
in Corollary~\ref{cor:mostly:monotone}.  Then, there exist
functions $f$, $g$, $q$, a corresponding product map $\hat \phi$,
and a subset $U_2 \subset B(R)$ with $\mu(U_2) > (1 -
256\delta^{1/4}) \mu(B(R))$ such that for $(x,y,z) \in U_2$,
\begin{displaymath}
d(\phi(x,y,z), \hat \phi(x,y,z)) = O(\epsilon R).
\end{displaymath}
\end{proposition}

\bold{Proof.} To simplify language, we assume that the dominant
orientation is right-side-up (the other case is identical). Let
$z_1$ (resp. $z_2$) denote the height of the bottom (resp. top) of
$B(R)$. If $(x,y,z) \in B(R)$, we let $\gamma_{xy}: [z_1,z_2] \to
B(R)$ denote the vertical geodesic segement $\gamma_{xy}(t) =
(x,y,t)$. Let $F_1$ (resp. $F_2$) denote the subset of the bottom
(resp. top) face of $B(R)$ which is within $8 \delta^{1/4}$ of a
point of $U$. Since $\mu(U) \ge (1-8\sqrt{\delta}) \mu(B(R))$,
each $F_i$ has nearly full $\mu$-measure. In fact if we let $U'
\subset B(R)$ denote the set of  points $(x,y,z)$ such that
$(x,y,z_1) \in F_1$, $(x,y,z_2) \in F_2$, and $\gamma_{xy}$ has
$\epsilon$-monotone image under $\phi$, then $\mu(U') \ge (1 - 8
\delta^{1/4})\mu(B(R))$. \mc{ say better and justify}

Note that $F_1$ is an $O(1)$ neighborhood of a (subset of a)
segment of a $x$-horocycle, say $\{ (x,y_1, z_1) \st x \in A \}$.
Since the restriction of $\phi$ to $U$ preserves the
$x$-horocycles, $\delta^{1/4} < \epsilon$, there exist numbers
$y_1'$ and $z_1'$  and a function $f: A \to \reals$ or $\ratls_m$
such that for $x \in A$, $\phi(x,y_1, z_1)$ is at most $O(\epsilon
R)$ distance from $(f(x),y_1', z_1')$. Similarly, $F_2$ is bounded
distance from a set of the form $\{ (x_2,y, z_2) \st y \in A' \}$,
and there exists a function $g: A' \to \reals$ or $\ratls_n$ such
that the restriction of $\phi$ to $F_2$ is $O(\epsilon R)$
distance from a map of the form $(x_2, y, z_2) \to
(x_2',g(y),z_2')$. \mc{ say better in
  DL case}

Let $U_1$ be as in Lemma~\ref{lemma:height:preserving}. Now
suppose $p = (x,y,z) \in U' \cap U_1$. Since $p \in U'$, $\phi(p)$
is $O(\epsilon R)$ from a vertical geodesic connecting a point in
the $O(\epsilon R)$ neighborhood of $\phi(x,y,z_1)$ to a point in
the $O(\epsilon R)$ neighborhood of $\phi(x,y,z_2)$. Hence,
$\phi(p)$ is within $O(\epsilon R)$ distance of the vertical
geodesic connecting $(f(x),y_1',z_1')$ to $(x_2',g(y),z_2')$.
This, combined with (\ref{eq:partial:product}) implies the
proposition, and hence Theorem~\ref{theorem:main:step:one}.
\qed\medskip

\noindent{\bf Remark:} The product map $\hat \phi$ produced in the
proof of Proposition \ref{prop:product:map} is not defined on the
entire box.  Since we are not assuming anything about the
regularity of the maps $f,g$ and $q$ which define $\hat \phi$, one
can choose an arbitrary extension to a product map defined on the
box. This is sufficient for our purposes here. \medskip

\noindent \bold{Order in which constants are chosen.}

\begin{itemize}
\item We may assume that $\epsilon$ is sufficiently small
  so that in Lemma~\ref{lemma:quadrilateral}, the
  $O(\epsilon r)$ error term is smaller then $(r/100)$.

\item We choose $N = N(\epsilon,\kappa,C)$, so that
  Lemma~\ref{lemma:subdivision} works. We may assume $N \in \zed$.

\item As described in \S\ref{sec:averaging}, we choose $\delta =
  \delta(\epsilon,\theta,\kappa,C)$, so that $\delta^{1/4} < \epsilon$
  (see proof of Proposition~\ref{prop:product:map}), $256\delta^{1/4}
  < \theta$ (see Proposition~\ref{prop:product:map}) and also
  $\delta^2 < \theta$, see Lemma~\ref{lemma:good:box}.

\item We choose $S = S(\delta, \kappa, \epsilon)$ so that $S >
  \frac{32 \kappa^3}{\epsilon \delta^4}$, as described in
  \S\ref{sec:averaging}.

\item For $s=1,\dots, S$, write $r_s = r_0 N^s$.

\item Write $L = \Delta r_0$. Choose
  $\Delta = N^{p}$ for some $p \in \zed$, and so that for $R = r_S$,
  the $O(R/L)$ error term in (\ref{eq:tiling:identity}) is at most
  $\delta^2$. Then the same is true for any $R = r_s$, $1 \le s \le S$.

\end{itemize}

Now assume $r_0$ is sufficiently large so that
Lemma~\ref{lemma:uniform:set} holds with $R = r_0$ and $\theta_1 =
2 \delta$. Theorem~\ref{theorem:main:step:one} holds with $\alpha =
1$ and $\beta = N^S$.

\section{Step II}
\label{sec:staircase}

In this section, we assume that $m > n$. We prove the following
theorem:
\begin{theorem}
\label{theorem:main:step:two} For every $\delta > 0$, $\kappa > 1$
and $C > 0$ there exists a constant $L_0 > 0$ (depending on
$\delta$, $\kappa$, $C$) such that the following holds: Suppose
$\phi: \X \to \pX$ is a $(\kappa,C)$ quasi-isometry. Then for
every $L > L_0$ and every box $B(L)$, there exists a subset $U
\subset B(L)$ with $|U| \ge (1 - \delta) |B(L)|$ and a
height-respecting map $\hat{\phi}(x,y,z) = (\psi(x,y,z),q(z))$
such that
\begin{itemize}
\item[(i)]
\begin{displaymath}
d(\phi|_U,  \hat{\phi}) = O(\delta L).
\end{displaymath}
\item[(ii)] For $z_1, z_2$ heights of two points in $B(L)$, we
have
\begin{equation}
\label{eq:q:bilishitz} \frac{1}{2\kappa} |z_1 - z_2| - O(\delta L)
< |q(z_1) - q(z_2)| \le 2 \kappa |z_1 - z_2| + O(\delta L).
\end{equation}
\item[(iii)] For all $x \in U$, at least $(1-\delta)$ fraction of
the vertical geodesics passing within $O(1)$ of $x$ are $(\eta,
O(\delta L))$-weakly monotone, where $\eta$ depends only on the model
space. 
\end{itemize}
\end{theorem}

\bold{Remark.} It is not difficult to conclude from
Theorem~\ref{theorem:main:step:two} that $\hat{\phi}$ is in fact a
product map (not merely height-respecting). However, we will not
need this.

Theorem~\ref{theorem:main:step:two} is true also for the case $m =
n$; its proof for that case is the content of \cite{EFW1}. The
proof presented in this section is much simpler, but applies only
to the case $m > n$.

 The main point of the proof is to show that if $m >
n$, then in the notation of Theorem~\ref{theorem:main:step:one}, for
each $i \in I_g$, the maps $\hat\phi_i$ must preserve the up
direction. This is done in \S\ref{sec:noflips}.  The deduction of
Theorem~\ref{theorem:main:step:two} from that fact is in
\S\ref{sec:lininguptheboxes}.  The argument here for showing that
the up direction is preserved uses the fact that when $m \neq n$,
each box $B(R)$ has most of it's mass at the bottom of the box.

\subsection{Volume estimates}
\label{sec:volume}

 This section collects a number of purely geometric
facts needed in the proof of Theorem \ref{theorem:main:step:two}.
The point is merely to show that quasi-isometries quasi-preserve
volume in an appropriate sense.

The following is a basic property shared for example by all
homogeneous spaces and all spaces with a transitive isometry group
(such as $\X$).

\begin{lemma}
\label{lemma:ball:growth} For $p \in \X$, let $D(p,r)$ denote the
metric ball of radius $r$ centered at $p$. Then, for every $b > a
> 0$ there exists $\omega = \omega(a,b)> 1$ with $\log \omega = O(b-a)$ such that for all $p,q \in \X$,
\begin{displaymath}
\omega(a,b)^{-1} |D(p,a)| \le |D(q, b)| \le \omega(a,b) |D(p,a)|,
\end{displaymath}
where $| \cdot |$ denotes volume (relative to the $\X$ metric).
Also $\log \omega(a,b) = O(b-a)$, where the implied constant
depends on the model space $\X$.
\end{lemma}

\bold{Proof.} The first statement is immediate since $\X$ is a
homogeneous space  and therefore $|D(q,a)|=|D(p,a)|$ for $p,q$. The second statement
is a consequence of exponential growth of balls.
\qed\medskip

In this section we prove some fairly elementary facts about
quasi-isometries quasi-preserving volume.  The main tool is the
following basic covering lemma.

\begin{lemma}
\label{lemma:covering} Let $X$ be a metric space and let $\mathcal
F$ be a collection of points in $X$.  Then for any $a>0$ there is
a subset $\mathcal G$ in $\mathcal F$ such that:
\begin{itemize}
\item[{\rm (i)}] The sets $\{B(x,a)|x{\in}\mathcal G\}$ are
pairwise
  disjoint.
\item[{\rm (ii)}] $\bigcup_{\mathcal
F}B(x,a)\subset\bigcup_{\mathcal G}B(x,5a).$
\end{itemize}
\end{lemma}

\noindent This lemma and it's proof (which consists of picking
$\cG$ by a greedy algorithm) can be found in \cite[Chapter 1]{He}.
This argument is implicit in almost any reference which discusses
covering lemmas.

Recall that we are assuming that $\phi$ is a continuous
$(\kappa,C)$ quasi-isometry.

\relax From this we can deduce the following fact about quasi-isometries
of $\X$.  This fact holds much more generally for metric measure
spaces which satisfy Lemma~\ref{lemma:ball:growth}.

\begin{proposition}
\label{prop:qivolume} Let $\phi:\X{\rightarrow}\pX$ be a
continuous $(\kappa,C)$ quasi-isometry.  Then for any $a \GG C$
there exists $\omega_1 > 1$ with $\log \omega_1 = O(a)$ such that
for any $U \subset \X$,
\begin{displaymath}
\omega_1^{-1} | \phi(N_a(U)) | \le | N_a(U)| \le \omega_1
|N_a(\phi(U))|
\end{displaymath}
where $N_a(U) = \{ x \in \X \st d(x,U) < a \}$.
\end{proposition}

\bold{Proof.} We assume that $a >  4 \kappa C$. Note that we are
assuming that every point is within distance $C$ of the image of
$\phi$. Let $\cF$ be the covering of $N_a(U)$ consisting of all
balls of radius $a$ centered in $U$. By
Lemma~\ref{lemma:covering}, we can find a (finite) subset $\cG$ of
$U$ such that $\bigcup_{x \in \cG}D(x,5a)$ cover $N_a(U)$ and such
that the balls centered at $\cG$ are pairwise disjoint. Hence,
\begin{displaymath}
\sum_{x \in \cG} |D(x,a)| \le | N_a(U) | \le \sum_{ x \in \cG}
|D(x,5a)|.
\end{displaymath}

Now $\phi(N_a(U))$ is covered by $\bigcup_{ x \in \cG} \phi(D(x,5
a)) \subset
  \bigcup_{ x \in \cG} D(\phi(x), 5 \kappa a + C)$.
Hence,
\begin{multline*}
|\phi(N_a(U))| \le \sum_{x \in \cG} | D(\phi(x), 5 \kappa a + C)|
\le \omega(a,5\kappa a +C ) \sum_{x \in \cG} | D(x,a)| \le \\ \le
\omega(a,5\kappa a + C) |N_a(U)|.
\end{multline*}
For the other inequality,
\begin{multline*}
|N_a(\phi(U))| \ge |N_a(\phi(\cG))| = \left| \bigcup_{x \in \cG}
D(\phi(x),a)\right| \ge \left| \bigcup_{x \in \cG}
D(\phi(x),a/\kappa - C)\right| \\ = \sum_{x \in \cG} |
D(\phi(x),a/\kappa - C)| \ge \omega(a/\kappa - 2 C, 5 a)^{-1}
\sum_{x \in G} | D(x,5 a)| \ge \omega_1^{-1} | N_a(U)|
\end{multline*}
\qed\medskip

\bold{Terminology.} The ``coarse volume''  of a set $E$ means the
volume of $N_a(E)$ for a suitable $a$. If the set $E$ is
essentially one dimensional (resp. two dimensional) we use the
term ``coarse length'' (resp. ``coarse area'') instead of coarse
volume but the meaning is still the volume of $N_a(E)$. We also
use $\ell(\cdot)$ to denote coarse length.

\subsection{The trapping lemma}

Again this section contains purely geometric facts
needed in the proof of Theorem \ref{theorem:main:step:two}.  The
facts in this section concern the geometry of the model space $\X$.
For a path $\gamma$, let
$\ell(\gamma)$ denote the length of $\gamma$ (measured in the
$\X$-metric).
Recall that we are assuming $m{\geq}n$.

\begin{lemma}
\label{lemma:length:area} Suppose $L$ is a constant $z$ plane,
and suppose $U$ is a bounded set contained in $L$. Suppose $k > r > 0$
and $\gamma$ is a path which stays at least $k$ units
below $L$, i.e. that $\max_{x{\in}\gamma}(h(x))<h(L)-k$.
Suppose also that any vertical
geodesic ray starting at $U$ and going down
intersects the $r$-neighborhood of
$\gamma$.  Then,
\begin{displaymath}
\ell(\gamma) \ge e^{c_1 k - c_2 r} \Area(U)
\end{displaymath}
where $c_1 > 0$ and $c_2 > 0$ depend only on the model space,
and both the length and the area are
measured using the $\X$ metric.

\end{lemma}

\mc{ This proof is written only for $DL(m,n)$.}

\bold{Proof.} We give a proof for $DL(m,n)$, the proof for $\Sol$
is similar. Let $\Delta$ denote the $r$-neighborhood of $\gamma$,
then $|\Delta| \le e^{c_2 r} \ell(\gamma)$. Pick $N$ so that
$\Delta$ stays above height $h(L)-N$. Let $\cA$ denote the set of
vertical segments of length $N$ which start at height $h(L)$ and
go down. Let $\cA_U$ denote the elements of $\cA$ which start at
points of $U$. Then $|\cA_U| = |U| m^N$. Now for $0<s<N$, any
point at height $h(U) - s$ intersects exactly $n^s m^{N-s}$
elements of $\cA$. Thus, by the assumption on the height of
$\Delta$, any point of $\Delta$ can intersect at most $n^k
m^{N-k}$ elements of $\cA$. But by assumption, $\Delta$ intersects
any element of $\cA_U$. Thus, $|\Delta| \ge (|U| m^N)/(n^k
m^{N-k}) = |U|(m/n)^k$, which implies the lemma. (Recall that in
our notation,  length = area = cardinality in $DL(m,n)$).

For $\Sol$ the proof is similar but uses smooth volume rather
than counting.  First observe that the volume of the $r$ neighborhood
of $\gamma$ is  at most $l(\gamma)e^{c_2r}$ where $c_2$ depends only on the model
space by exponential volume growth.
Second observe that projecting upward by $t$ units of height contracts
volume by $e^{c_1t}$.  Since $U$ is contained in the projection of
the $r$ neighborhood of $\gamma$ up to height $L$, and $\gamma$ is always
at least $k$ units below $L$, the desired estimate follows.

\qed\medskip

\bold{Remark:} When $m=n$, Lemma \ref{lemma:length:area} still holds (but
with $c_1 = 0$), and also in
addition with the word {\em below} replaced by the word {\em
above}. When $m{\neq}n$, volume decreases on upwards projection.

\subsection{Vertical Orientation preserved}
\label{sec:noflips}

Given a box $B(R)$, an $y$-horocycle $H$ in $B(R)$, and a number
$\rho$, we let the {\em shadow}, $\Sh(H,\rho)$, of $H$ in $B(R)$
be the set of points that can be reached by a vertical geodesics
going straight down from the $\rho$ neighborhood of $H$. Note that
if $H$ is the top of the box, $\Sh(H,1)$ is the entire box.
Similar definitions hold for $x$-horocycles, but then the shadow
will be above the horocycle.

The goal of this subsection is the following:
\begin{theorem}
\label{theorem:noflips} Suppose $m > n$ and that $\epsilon$ and
$\theta$ are sufficiently small (depending only on the model
space). Let $I$, $I_g$, $U_i$ and $\hat{\phi}_i$ be as in
Theorem~\ref{theorem:main:step:one}. Suppose $i \in I_g$. Then the
product map $\hat \phi_i :B(R) {\rightarrow} \pX$ can be written as
$\hat{\phi}_i(x,y,z) = (f_i(x),g_i(y),q_i(z))$, with $q_i: \reals \to
\reals$ coarsely orientation preserving.
\end{theorem}

\bold{Remark.} The result of Theorem~\ref{theorem:noflips} is {\em
  false} in the case $m = n$, since there exist ``flips'',
i.e. isometries which reverse vertical orientation. This is the
point where the proof in the case $m = n$ diverges from the case
$m > n$.
\medskip


In the rest of \S\ref{sec:noflips} we prove
Theorem~\ref{theorem:noflips}. We pick $i \in I_g$ and suppress
the index $i$ for the rest of this subsection.

Pick  $1 \GG \rho_2 \GG \rho_1 \GG \epsilon$ to be determined
later (see the end of this subsection).

\begin{lemma}
\label{lemma:goodshadow}  Let $\theta$ be as in Theorem \ref{theorem:main:step:one}. All but $O(4\sqrt{\theta})$ proportion of
the $y$-horocycles $H$ that are above the middle of the box $B(R)$
have all but $O(\sqrt{\theta})$ fraction of the $\mu$-measure of
both $N_{\rho_1 R}(H)$ and $\Sh(H,{\rho_1 R})$ in $U$.
\end{lemma}

\bold{Proof.} Let $P$ be a constant $z$ plane above the middle of
the box (and not too close to the top).  We choose horocycles
$H_i$ in $P$ such that $P=\coprod_{i}N_{\rho_1 R}(H_i)$. The
subset of $B(R)$ below $P$ is then the disjoint union of the
shadows $\coprod_i \Sh(H_i,\rho_1 R)$.  Since at least half the
measure of $B(R)$ is below $P$, it follows that there is some $i$
so that all but $O(\sqrt{\theta})$ of the $\mu$-measure of
$\Sh(H_i,\rho_1 R)$ is in $U$.  To guarantee the same fact about
$N_{\rho_1 R}(H)$, we pick $P$ such that $N_{\rho_1 R}(P)$ has all
but $O(\sqrt{\theta})$ fraction of its $\mu$-measure in $U$.
 \qed\medskip

\begin{lemma}
\label{lemma:goodplane} For any $H$ as in Lemma
\ref{lemma:goodshadow}, there exists a constant $z$ plane $P$ such
that $N_{\rho_1 R}(P) \cap \Sh(H,\rho_1 R){\cap} U$ contains all
but $O(\theta^{1/4})$ fraction of the $\mu$-measure in $N_{\rho_1
R}(P) \cap \Sh(H,\rho_1 R)$. Furthermore, we can choose $P$ and
$H$ such that $\rho_2 R < d(P,H) < 2 \rho_2 R$.
\end{lemma}

\bold{Proof.} Let
\begin{displaymath}
E = \Sh(H,\rho_1 R) \cap h^{-1}(h(H)- 2 \rho_2 R, h(H) - \rho_2
R).
\end{displaymath}
By Lemma~\ref{lemma:goodshadow}, $\mu(E \cap U) \ge (1 - c
\sqrt{\theta}/\rho_2) \mu(E) \ge (1 - \theta^{1/4}) \mu(E)$, where
$c$ is the implied constant in Lemma~\ref{lemma:goodshadow}, and
we have assumed that $c \theta^{1/2}/\rho_2 \le \theta^{1/4}$. Now
this is another application of Fubini's theorem, where we
partition $E$ into its intersections with neighborhood of
constant $z$ planes. \qed\medskip

\mc{One $H$ or almost all $H$?}

\begin{lemma}
\label{lemma:goodpieces} Let $P$ be as in the conclusion of Lemma
\ref{lemma:goodplane}.  There are subsets $S_1,S_2$ of
$P{\cap}B(R)$ such that
\begin{enumerate}
\item $N_{\rho_1 R}(S_i){\cap}U$ contains all but
$O(\theta^{1/4})$
  fraction of
the $\mu$-measure of $N_{\rho_1 R}(S_i)$. \item for $s_i{\in}S_i$
any path joining $s_1$ to $s_2$ of length less than $\kappa^3
\rho_2 R$ passes within $O(\rho_1 R)$ of $H$. \item For $i=1,2$,
$\Area(S_i) \GG \frac{1}{6}\Area(\Sh(H,\rho_1 R){\cap}P) > e^{c
\rho_2 R}\ell(H \cap B(R))$ where $c$ depends on the model spaces.
\end{enumerate}
\end{lemma}

\bold{Proof.} We divide $P{\cap}{\Sh(H,\rho_1 R)}$ into equal
thirds where each third has the entire $y$-extent and a third of
the $x$-extent. We let $\tilde{S}_1$ and $\tilde{S}_2$ be the two
non-middle thirds.  Now let $S_i$ be the portion of $\tilde{S}_i$
which is at least $\kappa^3 \rho_2 R$ away from the edges of
$B(R)$. The area of each of these regions is much more than the
coarse area $\ell(H \cap B(R))$ since projecting upwards decreases
area and each region projects upwards onto $H$. \qed\medskip

The proof of Theorem \ref{theorem:noflips} involves deriving
contradiction to the reversal of orientation of vertical geodesics
under $\hat \phi$ on $B_i(R)$. The goal is to show that if the
orientation were to reverse, we could find a path in the target
joining $\hat \phi(s_1)$ to $\hat \phi(s_2)$ that contradicts
Lemma \ref{lemma:goodpieces}.

\bold{Proof of Theorem \ref{theorem:noflips}.} We assume that
vertical orientation is not preserved but reversed. This means the
$z$ component $q(z)$ of the product map in
Theorem~\ref{theorem:main:step:one} is orientation reversing. Let
$H_p$ denote the $x$-horocycle through $p$. For $i=1,2$, let
\begin{displaymath}
S_i' = \{ p \in S_i \cap U \st \text{ for $j=1,2$, $\ell(H_p \cap
U
  \cap S_j) > 0.5 \ell(H_p \cap S_j)$ } \}
\end{displaymath}
Let $W_i = \hat{\phi}(S_i)$. By Proposition~\ref{prop:qivolume}, and
Lemma~\ref{lemma:goodpieces} part 3, we have
\begin{equation}
\label{eq:Wi:large:area} \Area(W_i) \ge e^{-D \rho_1 R}
\Area(S_i')
> e^{(c \rho_2 - D \rho_1)R} \ell(H)
\end{equation}
where $D$ depends only on the model spaces, $\kappa$ and $C$.
Assuming $\epsilon$ is sufficiently small, $\ell(H) \ll
\Area(W_i)$. Then by Lemma~\ref{lemma:length:area}, 99.9\% of the
geodesics going down from $W_i$ do not enter the $O(\rho_1 R)$
neighborhood of $\phi(H)$.  Let $W_i'$ denote the set of $p \in
W_i$ so that $99\%$ of the geodesics going down from (the
$1/2$-neighborhood of) $p$ do not enter the $O(\rho_1 R)$
neighborhood of $\phi(H)$. Then $\Area(W_i') \ge 0.9 \Area(W_i)$.

Let $\cH$ denote the set of $y$-horocycles $H'$ such that $W_i
\cap H'$ is non-empty for some (or equivalently for all) $i \in
\{1,2\}$. Let $\cH_i = \{H'\in H \st W_i'\cap H' \ne \emptyset
\}$. We claim that $\cH_1 \cap \cH_2 \ne \emptyset$. Indeed, $W_i
= \coprod_{H' \in H} W_i \cap H'$, and $\Area(W_i) = |\cH| c_i$,
where $c_i = |W_i \cap H'|$ is independent of $H' \in \cH$. Now,
\begin{displaymath}
0.9 c_i |\cH| = 0.9 \Area(W_i) \le \Area(W_i') \le c_i |\cH_i|.
\end{displaymath}
Thus, we have $|\cH_i| \ge 0.9 |\cH|$, and hence there exists $H'
\in \cH_1 \cap \cH_2$. By the definition of the $\cH_i$, we can
find $p_1$, $p_2$ such that $p_i \in H' \cap W_i'$. Then for
$i=1,2$, by the definition of $W_i'$, we can find geodesics
$\gamma_i$ going down from $p_i$, such that $\gamma_i$ do not
enter the $O(\rho_1 R)$ neighborhood  of
$\phi(H)$, and such that $\gamma_1$ and $\gamma_2$ meet at some
point $p'$. By construction $d(p', W_i) \le 2 \kappa^2 \rho_2 R$.
Concatenating subsegments of these two geodesics yields a path
connecting $p_1$ to $p_2$ of length $d(p_1,p_2)$ which avoids the
$O(\rho_1 R)$ neighborhood of $\phi(H)$.  Pulling back, we have a
path of length at most $16 \kappa^3 \rho_2 R$ and avoiding the
$O(\rho_1 R)$ neighborhood of $H$, which connects a point within
$O(\epsilon R)$ of $S_1$ to a point within
$O(\epsilon R)$ of $S_2$.  This contradicts Lemma
\ref{lemma:goodpieces}. \qed\medskip

\noindent{\bf Remark:} The only place in this paper where we make
essential use of the fact that $\phi$ satisfies $(2)$ of
Definition \ref{defn:qi} is in pulling back the path connecting
$p_1$ and $p_2$ at the end of the proof of Theorem
\ref{theorem:noflips}.

\bold{Choice of constants.} Let $A$ be the largest constant
depending only on $\kappa,C$ and the model spaces, which arises in
the course of the argument in \S\ref{sec:noflips}.

We choose $\rho_2$ so that $A \rho_2 < 1$. Similarly we choose
$\rho_1$ so that $A \rho_1 <\rho_2$ and $\epsilon$ $A \epsilon <
\rho_1$. We also choose $\theta$ so that $A \theta^{1/4} <
\rho_1$. We also make sure that $\epsilon$ and $\theta$ are
sufficiently small so that Theorem~\ref{theorem:noflips} applies.
In addition we choose $r_0$ in the statement of Theorem
\ref{theorem:main:step:one} such that the constant $e^{{(c \rho_2
- D\rho_1})R}$ that appears in the proof of Theorem
\ref{theorem:noflips} is at least $1000$.  Our other choices
guarantee that $\rho_2>\frac{D}{c}\rho_1$, so this is just a lower
bound on $R$ and therefore $r_0$.

\subsection{Proof of Theorem \ref{theorem:main:step:two}}
\label{sec:lininguptheboxes}

\bold{The uniform set and the exceptional set.} Let $I_g$ and
$U_i$, $i \in I_g$ be as in Theorem~\ref{theorem:main:step:one}. Let
$W = \bigcup_{i \in
  I_g} U_i$.

Recall that $Y_L$ is the set of vertical geodesics in $B(L)$. Here
we will work with a fixed geodesic $\gamma \in Y_L$.  Let $W^c
\subset B(L)$ denote the complement of $W$ in $B(L)$. For a point
$x \in \gamma$ and $T > 0$ let
\begin{displaymath}
P(x,\gamma, T) = | W^c \cap \gamma \cap D(x,T)|,
\end{displaymath}
where $D(x,T)$ is the ball of radius $T$ centered at $x$ (so that
$\gamma \cap D(x,T)$ is an interval of length $2T$ centered at
$x$).

\begin{lemma}
\label{lemma:uniform:implies:weakly:monotone} For every $\eta_1 >
0$ there exists $\eta > 0$ (with $\eta \to 0$ as $\eta_1 \to 0$)
such that the following holds: Suppose $\gamma$ is a geodesic ray
leaving $x$, and  for any $T > 1$, $P(x,\gamma, T)   < \eta_1 T$.
Then, $\bar{\gamma} = \phi \circ \gamma$ is
$(\eta,C_1)$-weakly-monotone, where $C_1 = O(\eta_1 R)$.
\end{lemma}

\bold{Proof.} Parametrize $\gamma$ so that $\gamma(0) = x$.
Without loss of generality, we may assume that $\gamma$ is going
up. Let $\bar{\gamma} = \phi \circ \gamma$. Suppose $0 < t_1 <
t_2$ are such that $h(\bar{\gamma}(t_1)) = h(\bar{\gamma}(t_2))$.
Write $q(t) = h(\bar{\gamma}(t)$. Subdivide $[t_1,t_2]$ into
intervals $I_1,\dots,I_N$ of length $\le \eta_1 R$ and so that the
length of all but the first and last is exactly $\eta_1 R$. We may
assume $N \ge 3$. \mc{ explain} Let $J \subset [1,\dots,N]$ be the
set of $j \in \zed$ such that $\bar{\gamma}(I_j) \cap W \ne
\emptyset$. For $j \in J$, pick $s_j$ such that $\bar{\gamma}(s_j)
\in W$, and pick $s\in I_j$ arbitrarily otherwise. Now
\begin{equation}
\label{eq:zero:sum} 0 = q(t_2) - q(t_1) = q(t_2) - q(s_{N'}) +
\sum_{\stackrel{j =3}{j \text{
      odd}}}^{N'} (q(s_j) - q(s_{j-2})) + q(s_1) - q(t_1),
\end{equation}
where $N'$ is either $N$ or $N-1$ depending on whether $N$ is odd
or even.

Let $Q_0 = \{ \text {odd } j \in [3,N'] \st \gamma(s_j) \in U_i, \
\gamma(s_{j-2}) \in U_i \}$ (same $U_i$). Let $Q_1$ denote the set
of odd $j \in [3,N']$ such that $\gamma(s_j)$ and
$\gamma(s_{j+1})$ are in different boxes $B_i(R)$. Finally, let
$Q_2$ denote the set of odd $j \in [3,N']$ such that $\gamma(I_j)
\subset W^c$ or $\gamma(I_{j-2}) \subset W^c$. By assumption,
$|Q_2| \le t_2/R$ and also $|Q_1| \le t_2/R$. \mc{ explain} Then,
$|Q| \ge (1/3) (t_2-t_1)/(\eta_1 R)$. Note that if $j \in Q$, then
$|q(s_j) - q(s_{j-2})| \ge \eta_1 R/(2 \kappa)$, and for any $j$,
$|q(s_{j+2}) - q(s_j)| \le 4 \kappa \eta_1 R$. Hence,
\begin{displaymath}
\sum_{j \in Q} q(s_j) - q(s_{j-2}) \ge |Q| \frac{\eta_1 R}{\kappa}
\ge \frac{t_2 - t_1}{6\kappa}.
\end{displaymath}
 Also,
\begin{displaymath}
\left| \sum_{j \in Q_1 \cup Q_2} q(s_{j}) - q(s_{j-2}) \right| \le
|Q_1 \cup Q_2| 2 \kappa \eta_1 R \le 2 \kappa \eta_1 t_2
\end{displaymath}
Plugging into (\ref{eq:zero:sum}) we see that
\begin{displaymath}
0 \ge \frac{t_2 - t_1}{6\kappa} - 2 \kappa \eta_1 t_2 - O(\eta_1
R),
\end{displaymath}
or
\begin{displaymath}
\frac{t_2 - t_1}{6\kappa} \le 2 \kappa \eta_1 t_2 + O(\eta_1 R)
\end{displaymath}
which implies the lemma.
\qed\medskip

Pick $A \GG 1$ (in fact we will eventually choose
$A=(4(128/\delta)^4$). Suppose $\gamma \in Y_L$, We define a point
$x \in \gamma$ to be {\em $A$-uniform along $\gamma$}, if for all
$T
> 1$,
\begin{displaymath}
\frac{P(x,\gamma,T)}{T} < A \frac{|\gamma \cap W^c|}{L}
\end{displaymath}

\color{black}
\begin{lemma}
\label{lemma:vitali} Let $\theta(\gamma)$ denote the proportion of
non-$A$-uniform points along $\gamma$. Then, $\theta(\gamma) \le
2/A$.
\end{lemma}

\bold{Proof.} This is a standard application of the Vitali
covering lemma. Let $\nu = \frac{|\gamma \cap W^c|}{L}$. Suppose
$x$ is non-uniform, then there is an interval $I_x$ centered at
$x$ such that
\begin{displaymath}
|I_x \cap W^c| \ge A \nu |I_x|.
\end{displaymath}
The intervals $I_x$ obviously cover the non-uniform set of
$\gamma$, and, by Vitali, we can choose a disjoint subset $I_j$
which cover at least half the measure of the non-uniform set.
Then,
\begin{displaymath}
| \bigcup I_j | \le  \sum |I_j| \le (A \nu)^{-1}|I_j \cap W^c| \le
(A \nu)^{-1} |\gamma \cap W^c|
\end{displaymath}
Dividing both sides by $L$ (the length $\gamma$), and recalling
that $|\gamma \cap W^c|/L = \nu$, we obtain the estimate.
\qed\medskip

Let $\theta_1 = \frac{\theta}{\eta_1} +\frac{2}{A}$.
\begin{corollary}
\label{cor:mostly:weakly:monotone} There exists a subset $U
\subset B(L)$ with $\mu(U) > (1 - 2 \sqrt{\theta_1}) \mu(B(L))$
such that for $x \in U$, $(1-\sqrt{\theta_1})$-fraction of the
geodesics passing within $(1/2)$ of $x$ are right-side-up
$(\eta,\eta_1 R)$-weakly-monotone.
\end{corollary}

\bold{Proof.} Let $Y'$ denote the space of pairs $(\gamma,x)$
where $\gamma \in Y_L$ is a vertical geodesic in $B(L)$ and $x \in
\gamma$ is a point. Let $E \subset Y'$ denote the set of pairs
$(\gamma,x)$ such that $|\gamma \cap W| \ge (1 - \eta_1/A) L$, and
$x$ is $A$-uniform along $\gamma$. Then, by
Lemma~\ref{lemma:vitali} we have $|E| \ge (1 - \theta_1) |Y'|$.
Let $U$ be the subset constructed by applying
Lemma~\ref{lemma:uniform:set}. Then $\mu(U) \ge (1 - 2
\sqrt{\theta_1})\mu(B(L))$, and for $x \in U$ by
Lemma~\ref{lemma:uniform:implies:weakly:monotone} at least
$(1-\sqrt{\theta_1})$ fraction of the geodesic rays leaving $x$
are $(\eta,O(\eta_1 R))$-weakly-monotone. \qed\medskip

\begin{lemma}
\label{lemma:weakly:monotone:to:height:preserving} Suppose $\phi$
and $B(L)$ and $U$ are as in
Corollary~\ref{cor:mostly:weakly:monotone}, and $\eta$ is
sufficiently small (depending only on the model space). Then,
there exist functions $\psi$, $q$, and a subset $U_1 \subset B(L)$
with $\mu(U_1) > (1 - 128\theta_1^{1/4}) \mu(B(L))$ such that for
$(x,y,z) \in U_1$,
\begin{equation}
\label{eq:partial:product:weak} d(\phi(x,y,z),(\psi(x,y,z),q(z))) =
O(\delta L)
\end{equation}
\end{lemma}

\bold{Proof.} This proof is identical to that of
Lemma~\ref{lemma:height:preserving}. \qed\medskip

\bold{Proof of Theorem~\ref{theorem:main:step:two}.} Choose $\eta$
so that Lemma~\ref{lemma:weakly:monotone:to:height:preserving}
holds. Choose $\eta_1$ so that
Lemma~\ref{lemma:uniform:implies:weakly:monotone} holds, and also
that the $O(\eta_1 R)$ term in
Lemma~\ref{lemma:uniform:implies:weakly:monotone} is at most
$\delta L$.  Choose $A^{-1} = (\delta/128)^4/4$ and choose $\theta =
(\delta/128)^4/\eta_1$ so that $128 \theta_1^4 < \delta$. Now the
theorem follows from combining
Corollary~\ref{cor:mostly:weakly:monotone} and
Lemma~\ref{lemma:weakly:monotone:to:height:preserving}.
\qed\medskip

\section{Step III}
\label{sec:fullmeasure}

In this section, we complete the proof of
Theorem~\ref{theorem:qisol} and Theorem~\ref{theorem:qidl}. We
assume that $\phi$ is a $\kappa,C$ quasi-isometry from $\X$ to
$\pX$ satisfying the conclusion of
Theorem~\ref{theorem:main:step:two}. All the arguments in this
section are valid also in the case $m=n$ (and are used in
\cite{EFW1}).

\subsection{A weak version of height preservation}
In this subsection, our main goal is to prove the following:

\begin{theorem}
\label{theorem:euc:plane:lin:drift} Let $\phi: \X \to \pX$ be a
$(\kappa,C)$ quasi-isometry satisfying the conclusions of
Theorem~\ref{theorem:main:step:two}. Then for any $\theta \ll 1$
there exists $M > 0$ (depending on $\theta,\kappa, C$) such that
for any $x$ and $y$ in $\X$ with $h(x) = h(y)$,
\begin{equation}
\label{eq:euc:plane:lin:drift} |h(\phi(x)) - h(\phi(y))| \le
\theta d(x,y) + M.
\end{equation}
\end{theorem}

\bold{Note.} This is a step forward, since the theorem asserts
that (\ref{eq:euc:plane:lin:drift}) holds {\em for all} pairs
$x,y$ of equal height (and not just on a set of large measure).

\medskip

We would like to restrict Theorem~\ref{theorem:main:step:two} to
the neighborhood of a constant $z$ plane. Let $\nu = \sqrt{\delta}$.
Fix a constant $z$ plane $P$. For notational convinience, assume that
$P$ is at height $0$. Let $R(L) \subset P$ denote the intersection
of $P$ with a box $B(2 L)$ whose top face is at height $L$ and
bottom face at $-L$. Then $R(L)$ is a rectangle (in fact, when
$m=n$, with this choice of $P$,  $R(L)$ is a square in the
euclidean metric). We will call $L$ the {\em size} of $R(L)$. Let
$R^+(L)$ denote the ``thickening'' of $R(L)$ in the $z$-direction
by the amount $\nu L$, i.e. $R^+(L)$ is the intersection of $B(2
L)$ with the region $\{p \in \X \st -(\nu/2)L \le h(p) \le
(\nu/2)L \}$, where as above $h(\cdot)$ denotes the height
function.

We now have the following corollary of
Theorem~\ref{theorem:main:step:two}:

\begin{corollary}
\label{cor:euc:plane:most:measure} Suppose $L > L_0$. Then for
every rectangle $R(L) \subset P$ there exists $U \subset R^+(L)$
with $\mu(U) \ge (1-\nu) \mu(R^+(L))$ and a standard map
$\hat{\phi}: U \to \pX$ such that $d(\phi|_U,\hat{\phi}) \le \nu
L$.  Furthermore, for any $p \in U$, for 99\% of the geodesics
$\gamma$ leaving $p$, $\phi(\gamma \cap B(2 L))$ is within $\delta
L$ of a vertical geodesic segment (in the right direction).
\end{corollary}

\bold{The tilings.} Choose $\beta \ll 1$ depending only on
$\kappa,C,m$ and $n$. When $m=n$, $\beta \approx
\frac{1}{\kappa^4}$. Let $L_j = (1+\beta)^j L_0$. For each $j
> 0$ we tile $P$ by rectangles $R$ of size $L_j$; we denote the
rectangles by $R_{j,k}$, $k \in \natls$. For $x \in \X$, let
$R_j[x]$ denote the unique rectangle in the $j$'th tiling to which
the orthogonal projection of $x$ to $P$ belongs.

\bold{Warning.} Despite the fact that $L_{j+1} = (1+\beta) L_j$,
the number of rectangles of the form $R_j[x]$ needed to cover a
rectangle of the form $R_{j+1}[y]$ is very large (on the order of
$e^{\beta L_j}$). This is because the Euclidean size of $R(L_j)$
is approximately $e^{L_j}$. \mc{ fix
  dependence on $m,n$}
\medskip

\bold{The sets $U_j$.} For each rectangle $R_{j,k}$,
Corollary~\ref{cor:euc:plane:most:measure} gives us a subset of
$R_{j,k}^+$ which we will denote by $U_{j,k}$. Let
\begin{displaymath}
U_j = \bigcup_{k=1}^\infty U_{j,k}.
\end{displaymath}
In view of Corollary~\ref{cor:euc:plane:most:measure}, for any $x
\in U_j$,
\begin{equation}
\label{eq:same:j} \sup_{y \in R_j^+[x] \cap U_j} |h(\phi(y)) -
h(\phi(x))| \le 2 \nu L_j
\end{equation}

We also have the following generalization:
\begin{lemma}
\label{lemma:next:level} For any $x \in U_j$ and any $y \in
R_{j+1}^+[x] \cap U_j$,
\begin{displaymath}
|h(\phi(y)) - h(\phi(x))| \le 12 \nu L_j.
\end{displaymath}
\end{lemma}

\bold{Proof.} Let $R_j[p]$ be a rectangle on the same ``row'' as
$R_j[x]$ and the same ``column'' as $R_j[y]$. Then, since $\nu \ll
1$, there exists an $x$-horocycle $H$ which intersects both
$R_j^+[x] \cap U_j$ and $R_j^+[p] \cap U_j$; let us denote the
points of intersection by $x_1$ and $p_1$ respectively.

Now for $i=1,2$ choose (sufficiently different) vertical geodesics
$\gamma_i$ coming down from (near) $x_1$ and $\gamma_i'$ coming
down from (near) $p_1$ such that for $i=1,2$, $\gamma_i(L_{j+1})$
and $\gamma_i'(L_{j+1})$ are close.  (here all the geodesics are
parametrized by arclength).  In view of
Corollary~\ref{cor:euc:plane:most:measure}, since $x_1$ and $p_1$
are in $U_j$, we may assume that there exist vertical geodesics
$\lambda_i$ and $\lambda_i'$ such that for $0 \le t \le L_j$,
$d(\gamma_i(t),\lambda_i) \le \nu+\eta t$ where $\eta \ll 1$.
Similarly, $d(\gamma_i'(t),\lambda_i') \le \nu+ \eta t$.

Thus, in particular, $h(\phi(\gamma_i(L_j))) \le h(\phi(x_1)) -
L_j/\kappa + \eta \le h(H) - L_j/(2 \kappa)$, and similarly
$h(\phi(\gamma_i'(L_j))) \le h(p_1) - L_j/(2 \kappa)$. Now note
that $d(\gamma_i(L_j),\gamma_i'(L_j) = \beta L_j +O(1)$.
Hence $d(\phi(\gamma_i(L_j)),\phi(\gamma_i'(L_j))) \le 2 \kappa \beta
L_j + O(1)$, and by assumption $\kappa^2 \beta \ll 1$.
Then by
Lemma~\ref{lemma:quadrilateral}, $\phi(x_1)$ and $\phi(p_1)$ are
near the same horocycle, and thus, in particular,
\begin{equation}
\label{eq:x:hor} |h(\phi(x_1)) - h(\phi(p_1))| \le 4 \nu L_j
\end{equation}
Similarly, we can find $p_2 \in R_j^+[p] \cap U_j$ and $y_2 \in
R_j^+[y] \cap U_j$ such that $p_2$ and $y_2$ are on the same
$y$-horocycle. Then, by the same argument,
\begin{equation}
\label{eq:y:hor} |h(\phi(p_2)) - h(\phi(y_2))| \le 4 \nu L_j
\end{equation}
Hence, in view of (\ref{eq:x:hor}), (\ref{eq:y:hor}), and
(\ref{eq:same:j}),
\begin{displaymath}
|h(\phi(x)) - h(\phi(y))| \le 12 \nu L_j,
\end{displaymath}
as required. \qed\medskip

\begin{lemma}
\label{lemma:Uj:Ujplus:one} Suppose $p \in R_j^+[x] \cap U_j$, $q
\in R_{j+1}^+[x] \cap U_{j+1}$. Then,
\begin{equation}
\label{eq:Uj:Ujplus:one} |h(\phi(p)) - h(\phi(q))| \le 16 \nu L_{j+1}
\end{equation}
\end{lemma}

\bold{Proof.} Note that the orthogonal projection of $U_j \cap
R_{j+1}^+[x]$ to $R_{j+1}[x]$ has full $\mu$-measure (up to order
$\nu$). The same is true of $U_{j+1} \cap R_{j+1}^+[x]$. Thus, the
projections intersect, and thus we can find $p' \in U_j \cap
R_{j+1}^+[x]$ and $q' \in R_{j+1}^+[x] \cap U_{j+1}$ such that
$d(p',q') \le 2 \nu L_{j+1}$. Now, in view of
Lemma~\ref{lemma:next:level},
\begin{displaymath}
|h(\phi(p)) - h(\phi(p'))| \le 12 \nu L_j
\end{displaymath}
and in view of (\ref{eq:same:j}),
\begin{displaymath}
|h(\phi(q')) - h(\phi(q))| \le 2 \nu L_{j+1}
\end{displaymath}
This implies (\ref{eq:Uj:Ujplus:one}). \qed\medskip

\bold{Proof of Theorem~\ref{theorem:euc:plane:lin:drift}.} We have
\begin{displaymath}
R_0[x] \subset R_1[x] \subset R_2[x] \subset \dots
\end{displaymath}
and
\begin{displaymath}
R_0[y] \subset R_1[y] \subset R_2[y] \subset \dots
\end{displaymath}
There exists $N$ with $L_N$ comparable to $d(x,y)$ such that
(after possibly shifting the $N$'th grid by a bit) $R_N[x] =
R_N[y]$. Now for $0 \le j \le N$, pick $x_j \in R_j^+[x] \cap
U_j$, $y_j \in R_j^+[y] \cap U_j$. We may assume that $x_N = y_N$.
Now, using Lemma~\ref{lemma:Uj:Ujplus:one},
\begin{align*}
|h(\phi(x_0)) - h(\phi(y_0))| & \le \sum_{j=0}^{N-1}
|h(\phi(x_{j+1})-h(\phi(x_j))| + \sum_{j=0}^{N-1}
|h(\phi(y_{j+1})-h(\phi(y_j))|   \\
& \le 2 \sum_{j=0}^{N-1} 16 \nu L_{j+1} \\
& \le \frac{32 \nu}{\beta} L_{N},
\end{align*}
where in the last line we used that $L_j = (1+\beta)^j L_0$. Now
since $x_0 \in R_0[x]$, $d(x,x_0) \le L_0$, so $|h(\phi(x) -
\phi(x_0)| = O(L_0)$. Similarly, $|h(\phi(y) - h(\phi(y_0))| =
O(L_0)$. Also note that $L_{N+1}$ is within a factor of $2$ of
$d(x,y)$. Thus the theorem follows. \qed\medskip





\subsection{Completion of the proof of height preservation}

\begin{lemma}
\label{lemma:all:weakly:monotone} Let $\phi: \X \to \pX$ be a
$(\kappa,C)$ quasi-isometry. Then for any $\eta \ll 1$ there
exists $C_1 > 0$ (depending on $\eta, \kappa, C$) such that for
any vertical geodesic ray $\gamma$, $\phi \circ \gamma$ is $(\eta,
C_1)$-weakly monotone.
\end{lemma}

\bold{Proof.} This is a corollary of
Theorem~\ref{theorem:euc:plane:lin:drift}. Suppose $\gamma$ is a
vertical geodesic ray parametrized by arclength, and $\bar{\gamma}
= \phi \circ \gamma$. Suppose $0 < t_1 < t_2$ are such that
$h(\bar{\gamma}(t_1)) = h(\bar{\gamma}(t_2))$. We now apply
Theorem~\ref{theorem:euc:plane:lin:drift} to $\phi^{-1}$ instead
of $\phi$ (with $x = \bar{\gamma}(t_1))$ and $y =
\bar{\gamma}(t_2))$. We get $|h(\gamma(t_1))-h(\gamma(t_2))| \le
\theta d(\bar{\gamma}(t_1)), \bar{\gamma}(t_2)) + O(M)$, i.e.
\begin{displaymath}
|t_2 - t_1| \le \theta \kappa^2 |t_2 - t_1| + O(M)
\end{displaymath}
I.e. $\bar \gamma$ is $(\theta\kappa^2, O(M))$-weakly monotone.
\qed\medskip


\bold{Proof of Theorem~\ref{theorem:qisol} and
  Theorem~\ref{theorem:qidl}.}
Suppose $p_1$ and $p_2$ are two points of $\X$, with $h(p_1) =
h(p_2)$. We can find $q_1$, $q_2$ in $\X$ such that
$p_1,p_2,q_1,q_2$ form a quadrilateral. By
Lemma~\ref{lemma:all:weakly:monotone}, each of the segments
$\gamma_{ij}$ connecting a point in the $O(1)$ neighborhood of
$p_i$ to a point in the $O(1)$ neighborhood of $q_j$ maps under
$\phi$ to an $O(\eta,C_1)$-weakly monotone quasi-geodesic segment.
Then by Lemma~\ref{lemma:epsilon:monotone:segment}, and
Lemma~\ref{lemma:quadrilateral}, we see that $h(\phi(p_1)) =
h(\phi(p_2)) + O(C_1)$. \qed\medskip

\section{Deduction of rigidity results}
\label{section:atinfinity}

The purpose of this section is to apply the previous results on
self quasi-isometries of $\Sol$ and the $\DL$-graphs to understand
all finitely generated groups  quasi-isometric to either one.
This follows a standard outline: if $\G$ is quasi-isometric to $X$
then $\G$ quasi-acts on $X$ (in this case that just means there is
a homomorphism $\G \to \QI(X)$ with uniformly bounded constants).
We then need to show that such a quasi-action can be conjugated to
an isometric action.   The basic ingredients we need to do this
are the following:

\begin{theorem}\cite{FM2}\label{qsr} Every uniform quasi-similarity action on $\reals$ is
bilipschitz conjugate to a similarity action.
\end{theorem}

\noindent The proof of this theorem makes substantial use of work
of Hinkannen \cite{H} who had shown that a uniform quasi-symmetric
action was quasi-symmetrically conjugate to a symmetric action.

\begin{theorem}\cite{MSW}\label{qsqp} Let $\G$ have a uniform quasi-similarity action on
$\ratls_m$.   If the $\G$ action is cocompact on the space of pairs of distinct points in
$\ratls_m$ then there is some $n$ and a similarity action of $\G$ on $\ratls_n$ which is
bilipschitz conjugate to the given quasi-similarity action.
\end{theorem}

It is useful to think about these results in a quasi-action interpretation.   One can
view $\reals$ as $S^1 - \{pt\}$, and interpret a uniform quasi-similarity action on
$\reals$ as the boundary of a quasi-action on $\half^2$ fixing a point at infinity.   The
result of Farb and Mosher then says that this quasi-action is quasi-conjugate to an
isometric action on $\half^2$.    The interpretation of the second result is similar,
with a tree of valence $m+1$ replacing $\half^2$.     The hypothesis of cocompactness on
pairs in that theorem then translates to cocompactness of the quasi-action on the tree.

We now state and prove a result that immediately implies Theorem
\ref{theorem:nolattice}.  This result is also used in \cite{EFW1}.

\begin{theorem}
\label{theorem:qitoisom} Assume every $(\kappa,C)$ self
quasi-isometry of $\Sol$ is at bounded distance from a
$b$-standard map where $b=b(\kappa,C)$. Then any uniform group of
quasi-isometries of $\Sol$ is virtually a lattice in $\Sol$.
\end{theorem}

\bold{Proof.} Let $f: \G \to \Sol$ be a quasi-isometry.   For each
$\g$ in $\G$ we have the self-quasi-isometry $T_\g$ of $\Sol$
given by $$ x \mapsto f(\g f^{-1}(x)) $$

\noindent By Theorem \ref{theorem:qisol}, $T_\g$ is bounded
distance from a standard map.  On a subgroup $\G'$ of $\G$ of
index at most two, this gives a homomorphism $\Phi: \G' \to
\QSim(\reals) \times \QSim(\reals)$.    By Theorem \ref{qsr}, each
of these quasi-similarity actions on $\reals$ can be bilipschitz
conjugated to a similarity action.   This gives $\Psi: \G' \to
\Sim(\reals) \times \Sim(\reals)$.

Since the quasi-isometries $T_\g$ have uniformly bounded constants, we know that the
stretch factors of the two quasi-similarity actions $\Phi$ are approximately on the curve
$(e^{mt},e^{-nt})$ - meaning that the products weighted by these factors are uniformly
close to $1$.    This therefore is true for $\Psi$ as well.   So, in the sequence:

$$ \G' \to Sim(\reals)\times Sim(\reals) \to \reals \times \reals$$

\noindent where the final map is the log of the stretch factor, we
know that the image lies within a bounded neighborhood of the line
$ny = -mx$.   Since the image is a subgroup, this implies it must
lie on this line.  Since the subgroup of $\Sim(\reals) \times
\Sim(\reals)$ above this line is $\Sol$, we have produced a
homorphism

$$ \Psi : \G' \to \Sol$$

\noindent We now show that the kernel is finite and the image
discrete and cocompact.  This follows essentially from the fact
that the map $f$ is a quasi-isometry.

Consider a compact subset $K \subset  \Sol$.    The set $F = \Psi^{-1}(K)$ consists of
maps with uniformly bounded stretch factors, and which move the origin at most a bounded
amount.   Transporting this information back to the standard maps of $\Sol$, we see that
for $\g \in F$ the maps $T_\g$ move the identity a uniformly bounded amount.    However,
the quasi-action $T$ of $\G$ on $\Sol$ is the $f$-conjugate of the left action of $\G$ on
$\G$.  This action is proper, so we conclude that $F$ is finite.   This implies that
$\Psi$ has finite kernel and discrete image.   In the same way, the fact that the $\G$
action on $\G$ is transitive implies that the image of $\Psi$ is cocompact.

Thus the image of $\G'$ is a lattice in $\Sol$.

\qed {\medskip}

\noindent This proves Theorem \ref{theorem:nolattice}, since if $m
\neq n$, the group $\Sol$ is not unimodular and therefore does not
contain lattices.

We next prove Theorem \ref{theorem:dl}. In fact, we show

\begin{theorem}
\label{theorem:qitoisomDL} Assume every $(\kappa,C)$ self
quasi-isometry of $\DL(m,n)$ is at bounded distance from a
$b$-standard map where $b=b(\kappa,C)$. Then any uniform group of
quasi-isometries of $\DL(m,n)$ is virtually a lattice in
$\Isom(\DL(n',n'))$ where $n',m,n$ are all powers of a common
integer.
\end{theorem}

Some complications arise from the differences between Theorem
\ref{qsqp} and Theorem \ref {qsr}. We need the following theorem
of Cooper:

\begin{theorem}\cite{Co}
\label{theorem:cooper} The metric spaces $\ratls_p$ and $\ratls_q$
are bilipschitz equivalent if and only if there are integers
$d,s,t$ so that $p=d^s$ and $q=d^t$.
\end{theorem}

This immediately implies a weaker version Theorem
\ref{theorem:different:dl:not:qi}. We now turn to theorem
\ref{theorem:qitoisomDL}.

\bold{Proof.} We proceed as in the previous proof for $\Sol$. The
first difference is that to apply theorem \ref{qsqp} we need to
know that the quasi-similarity actions of $\G'$ on $\ratls_n$ and
$\ratls_m$ are cocompact on pairs of points.   As discussed above,
this is equivalent to asking the corresponding quasi-action on the
trees of valence $n+1$ and $m+1$ to be cocompact.    This then
follows immediately from the fact that $\G'$ is cocompact on
$\DL(m,n)$.

Thus we have $\Psi: \G' \to \Sim(\ratls_a) \times \Sim(\ratls_b)$
for some $a$ and $b$. Thus we know that we have $d_i,s_i,t_i$ for
$i=1,2$ with $n={d_1}^{s_1}$, $m={d_2}^{s_2}$ and:

$$\Psi : \G' \to \Sim(\ratls_{{d_1}^{t_1}}) \times \Sim(\ratls_{{d_2}^{t_2}})$$

\noindent We know, as before, that the weighted stretch factors
are approximate inverses.  In this case the stretch factors are in
$\integers$ - in $Sim(\ratls_m)$ one can stretch only by powers of
$m$.   Thus the image is a subgroup lying on the line $\{(a,b);
a*\log{d_1}*\frac{t_1}{s_1} + b*\log{d_2}*\frac{t_2}{s_2} =0\}$.
For this to be a non-empty subgroup  of $\integers^2$ we must have
$\frac{\log{d_1}}{\log{d_2}}$ rational, which implies that  there
is a $d$ with $d_1=d^u$, $d_2=d^v$ for some $u$ and $v$.

There is still some ambiguity in the choices, since many groups
occur as subgroups of $\Sim(\ratls_{p^k})$ for many different $k$.
As in the construction of \cite{MSW} we can make the choices
unique by choosing the $t_i$ the maximum possible, so that all
powers of ${d_i }^ {t_i}$ occur as stretch factors. With these
choices we are forced to have the line $\{(a,b): a+b=0\}$ as this
is the only line of negative slope in $\integers^2$ surjecting to
both factors.  Thus we have $\Psi : \G' \to \Sim(\ratls_{d^{t_1}})
\times \Sim(\ratls_{d^{t_2}})$, with the image contained in the
subgroup having inverse stretch factors. This group is, up to
finite index, $\Isom(\DL(d^{t_1},d^{t_2})$. So we have:

$$\Psi : \G' \to \Isom(\DL(d^{t_1},d^{t_2})$$

\noindent Exactly as before, one can see that the kernel is finite
and the image is a lattice, which implies that $t_1=t_2$.   This
implies that $\G$ is amenable, and hence it and $\DL(m,n)$ have
metric F\"olner sets.      This is true only for $m=n$, which
completes the proof. \qed\medskip

\noindent This immediately implies Theorem \ref{theorem:dl}, since
$\DL(m,n)$ is only amenable as a metric space when $m=n$.

\bold{Proof of Theorem~\ref{theorem:different:sol:not:qi}.} Since
all $\Solnn$ are obviously quasi-isometric to one another, it
suffices to consider the case $m{\neq}n$. This then follows
immediately from Theorem~\ref{theorem:qisol} and \cite[Theorem
5.1]{FM3}. \qed\medskip

\begin{proposition}
Theorem \ref{theorem:qidl} implies Theorem
\ref{theorem:different:dl:not:qi}.
\end{proposition}

\bold{Proof.} In view of Theorem~\ref{theorem:qidl}, the proof of
this result is similar to the last one.  The point is that (up to
permuting $m$ and $n$) the quasi-isometry
$\DL(m,n){\rightarrow}\DL(m',n')$ induces quasi-similarities
$\ratls_n{\rightarrow} \ratls_{n'}$ and
$\ratls_m{\rightarrow}\ratls_{m'}$.  Theorem \ref{theorem:cooper}
then implies that $m$ and $m'$ are both powers of some number $d$
and that $n$ and $n'$ are both powers of some number $s$. However,
since the quasi-similarities both come from the same map on
vertical geodesics, the scale factors must agree.  This
immediately implies $\log m'/\log m = \log n'/\log n$.
\qed

\noindent  Department of Mathematics, University of Chicago, Eckhart
Hall, 5734 S. University Avenue, Chicago, Illinois 60637.

\medskip
\noindent Department of Mathematics, Indiana University, Rawles
Hall, Bloomington, IN, 47405.

\medskip
\noindent Department of Mathematics, Statistics, \& Computer
Science, University of Illinois at Chicago\& 322 Science \&
Engineering Offices (M/C 249), 851 S. Morgan Street Chicago, IL
60607-7045.\

\end{document}